\documentclass[journal]{IEEEtran}

\usepackage{etex}

\usepackage{cite}
\usepackage{amsmath}
\usepackage{amssymb}
\usepackage{balance}
\usepackage{graphicx}
\usepackage{subfigure}
\usepackage{epstopdf}
\usepackage{amsfonts}
\usepackage[all]{xy}
\usepackage{array,xcolor}
\usepackage{color}
\usepackage{bbm}
\usepackage{dsfont}
\usepackage{caption}
\usepackage{booktabs}
\usepackage{multirow}

\newcommand*{\qed}{\hfill\ensuremath{\blacksquare}}%

\usepackage[font=small,labelfont=bf]{caption}

\usepackage{tikz}
\usetikzlibrary{shapes,matrix,decorations.shapes}
\usetikzlibrary{arrows,decorations.pathmorphing,backgrounds,positioning,fit,matrix}
\usetikzlibrary{shapes,arrows,chains}
\usetikzlibrary[calc]
\usetikzlibrary[plotmarks]
\usetikzlibrary[patterns,decorations]

\usepackage{algorithm}
\usepackage{algpseudocode}
\algrenewcommand\algorithmicrequire{\textbf{Input:}}
\algrenewcommand\algorithmicensure{\textbf{Output:}}

\usepackage[colorlinks=true,      
			linkcolor=black,      
			citecolor=black,      
			filecolor=black,      
			urlcolor=blue ]{hyperref}

\newcolumntype{C}[1]{>{\centering}m{#1}}
\newtheorem{definition}{Definition}

\newtheorem{proposition}{Proposition}
\newtheorem{remark}{Remark}
\newtheorem{lemma}{Lemma}
\newtheorem{example}{Example}
\newtheorem{problem}{Problem}

\newtheorem{assumption}{Assumption}

\newcommand{\tr}{{{\mathsf T}}}

\newcolumntype{C}[1]{>{\centering}m{#1}}

\newcounter{tempEquationCounter}
\newcounter{thisEquationNumber}

\begin{document}

\title{Distributed Design for Decentralized Control using Chordal Decomposition and ADMM}
\author{Yang~Zheng, Maryam Kamgarpour, Aivar Sootla, and Antonis~Papachristodoulou
\thanks{Y. Zheng is supported in part by the Clarendon Scholarship, and in part by the Jason Hu Scholarship. A. Sootla and A. Papachristodoulou are supported by EPSRC Grant EP/M002454/1. The work of M. Kamgarpour is gratefully supported by ERC Starting Grant CONENE.}
\thanks{Y. Zheng was with the Department of Engineering Science at the University of Oxford. He is now with the SEAS and CGBC at Harvard University. Email: zhengy@g.harvard.edu.}
\thanks{M. Kamgarpour is with the Department of Electrical Engineering and Information Technology at ETH Zurich. Email: mkamgar@control.ee.ethz.ch.}
\thanks{A. Sootla and A. Papachristodoulou are with the Department of Engineering Science at the University of Oxford. Emails: \{aivar.sootla, antonis\}@eng.ox.ac.uk.}
}

\maketitle

\begin{abstract}
    We propose a distributed design method for decentralized control by exploiting the underlying sparsity properties of the problem.  Our method is based on chordal decomposition of sparse block matrices and the alternating direction method of multipliers (ADMM). We first apply a classical parameterization technique to restrict the optimal decentralized control into a convex problem that inherits the sparsity pattern of the original problem. The parameterization relies on a notion of strongly decentralized stabilization, and sufficient conditions are discussed to guarantee this notion. Then, chordal decomposition allows us to decompose the convex restriction into a problem with partially coupled constraints, and the framework of ADMM enables us to solve the decomposed problem in a distributed fashion. Consequently, the subsystems only need to share their model data with their direct neighbours, not needing a central computation.     Numerical experiments demonstrate the effectiveness of the proposed method.
\end{abstract}


\IEEEpeerreviewmaketitle

\section{Introduction}

Many real-world complex systems, such as aircraft formation, automated highways and power systems, consist of a large number of interconnected subsystems. {Often in these interconnected systems, the controllers have only access to each subsystem's state information. The problem of design of stabilizing and optimal controllers based on only each subsystem's state information is referred to as \emph{decentralized control}. Due to its importance, this problem has attracted research attention since the late seventies~\cite{siljak2011decentralized,lunze1992feedback}.}

Early efforts have centered on decentralized stabilization and its algebraic characterization through the concept of \emph{decentralized fixed modes}~\cite{wang1973stabilization}. {These are the set of eigenvalues that remain unchanged under any decentralized feedback.} One seminal result is that a system is stabilizable by a decentralized controller if and only if its decentralized fixed modes have negative real parts~\cite{wang1973stabilization}. Since then, a wide range extensions of decentralized control have been investigated, either by considering various types of performance guarantees {in addition to stability}~\cite{geromel1994decentralized}, or by taking into account neighbouring information for feedback, known as distributed control~\cite{jovanovic2016controller}. Several classes of systems have been identified, which allow convex formulations for the design of distributed $\mathcal{H}_{\infty}$ and $\mathcal{H}_2$ controllers, \emph{e.g.}, quadratically invariant systems~\cite{rotkowitz2006characterization}. Also, some numerical approaches have been proposed to find an approximate solution to the optimal distributed control problem~\cite{lin2011augmented,fazelnia2017convex}. The case of decentralized control in the presentence of input and state constraints is addressed in~\cite{furieri2017value}.

A common assumption made in these papers is that a central model of the global plant is available, indicating that the design is performed in a \emph{centralized} fashion even though the implementation of controllers is \emph{decentralized}. However, this may be impractical for certain complex systems that are shared between private individuals, such as transportation systems and power-grids. In this case, a complete model may not be available due to privacy concerns of model information for the subsystems.   %
While discussions on distributed design relying on limited model information can be traced back to~\cite[Chapter 9]{lunze1992feedback}, practical approaches to this problem are an active topic of current research. For example, performance bounds of designing linear quadratic regulators distributedly were discussed for systems with invertible input matrix in~\cite{langbort2010distributed}.  The distributed design framework of~\cite{langbort2010distributed} has been used to discuss the best closed-loop performance achievable by distributed design strategies for a class of fully actuated discrete-time systems~\cite{farokhi2013optimal}. In~\cite{shah2013optimal}, independent decoupled problems were derived for optimal decentralized control by utilizing the properties of posets. Recent work has started to use distributed optimization techniques to realize distributed synthesis in the dissipative framework~\cite{ahmadi2018distributed}. Furthermore, the recently proposed system level approach has been promoted to address distributed design of dynamic distributed controllers~\cite{wang2017separable}.

In this paper, we propose a new distributed design method for optimal decentralized control by exploiting the sparsity structure of the system. {Our method uses local information on system model to design controllers that rely on subsystems' state measurements.} The idea originates from the connection between sparse positive semidefinite (PSD) matrices and chordal graphs~\cite{agler1988positive, grone1984positive}. The celebrated chordal decomposition in graph theory~\cite{agler1988positive, grone1984positive} allows us to decompose a large sparse PSD cone into a set of smaller and coupled ones, and has been successfully applied to decompose sparse semidefinite programs (SDPs)~\cite{fukuda2001exploiting,andersen2010implementation}. These results have recently been used for performance analysis of sparse linear systems~\cite{mason2014chordal,ZKSP2017scalable,andersen2014robust}, leading to significantly faster solutions than using standard dense methods. Despite scalability of these approaches, they all required global model information for centralized computation. The authors in~\cite{Zheng2017Scalable} proposed a sequential approach to improve the scalability of solving a stabilization problem of networked systems, where model privacy can be maintained as a byproduct.

This paper extends the scope of exploiting chordal decomposition to distributed design of optimal decentralized control. By using a classical parameterization technique that relies on a notion of \emph{strongly decentralized stabilization}~\cite{geromel1994decentralized}, the optimal decentralized control can be restricted to a convex problem that inherits the original sparsity pattern in the system. The convex restriction can be equivalently decomposed into a problem with partially coupled constraints, and we introduce a distributed algorithm to solve the decomposed problem based on the framework of \emph{alternative direction method of multipliers} (ADMM). Precisely, the main contributions of this paper are:
\begin{itemize}
  \item {We provide sufficient conditions to guarantee the feasibility of the proposed convex restriction. These conditions are based on characterizing the cases in which   the closed-loop system with decentralized feedback admit a block-diagonal Lyapunov function\footnote{The authors have summarized some preliminary results in an unpublished technical report~\cite[Section 3]{zheng2017convex}. The current manuscript serves as the official version of the report~\cite{zheng2017convex}, and we do not consider~\cite{zheng2017convex} for publication.}. In particular, we identify two classes of networked systems   admitting strongly decentralized stabilization.}
  \item One notable feature of the convex restriction is that the original sparsity pattern of the system is inherited in the resulting convex optimization problem. We combine chordal decomposition with ADMM to solve the convex problem in a distributed fashion. In our algorithm, no central model of the global plant is required and the subsystems only need to share their model data with their neighbours, which help preserve the privacy of model data.
\end{itemize}

The rest of this paper is organized as follows. We present the problem formulation in Section~\ref{Section:ProblemStatement}. In Section~\ref{Section:SufficientCondition}, we discuss sufficient conditions on strongly decentralized stabilization. Section~\ref{Section:Decomposition} applies a chordal decomposition technique to derive a decomposed problem, and a distributed algorithm is introduced to solve the decomposed problem in Section~\ref{Section:ADMM}. Numerical examples are given in Section~\ref{Section:Results}. We conclude this paper in Section~\ref{Section:Conclusion}.

\section{Background and Problem Statement}
\label{Section:ProblemStatement}

\subsection{Optimal decentralized control}

A directed graph $\mathcal{G}(\mathcal{V},\mathcal{E})$ is defined by a set of nodes $\mathcal{V}=\{1,2,\dots, N\}$ and a set of edges $\mathcal{E} \subseteq \mathcal{V} \times \mathcal{V}$. We consider a complex system consisting of $N$ subsystems. The interactions between subsystems are modeled by a plant graph $\mathcal{G}_p(\mathcal{V},\mathcal{E}_p)$, in which each node in $\mathcal{V}$ denotes a subsystem, and the edge $(i,j) \in \mathcal{E}_p$ means that subsystem $i$ has dynamical influence on subsystem $j$.
The dynamics of each subsystem $i \in \mathcal{V}$ are
\begin{equation}\label{E:SubDynamics}
    \dot{x}_i(t) = A_{ii}x_i(t) + \sum_{j \in \mathbb{N}_i} A_{ij}x_j(t) + B_i u_i(t) + M_i d_i(t),
\end{equation}
where $x_i \in \mathbb{R}^{n_i}, u_i \in \mathbb{R}^{m_i}, d_i \in \mathbb{R}^{q_i}$ denote the local state, input and disturbance of subsystem $i$, respectively, and $\mathbb{N}_i$ denotes the set of neighbouring nodes that influence node $i$, \emph{i.e.},
$
    \mathbb{N}_i = \{j \mid (j,i) \in \mathcal{E}\}.
$
In~\eqref{E:SubDynamics}, $A_{ii} \in \mathbb{R}^{n_i \times n_i}, B_i \in \mathbb{R}^{n_i \times m_i}, M_i \in \mathbb{R}^{n_i \times q_i}$ represent local dynamics, and $A_{ij} \in \mathbb{R}^{n_i \times n_j}$ represents the interaction with neighbors. In this paper, we refer to $A_{ii}, B_i, M_i, A_{ij}$ as \emph{model data} of the system.

By collecting the subsystems' states, the overall system can be described compactly as
\begin{equation}\label{E:GlobalDynamics}
    \dot{x}(t) = Ax(t) + Bu(t) + Md(t),
\end{equation}
where $x:=[x_1^{\tr},x_2^{\tr}, \ldots, x_N^{\tr}]^{\tr}$, and the vectors $u, d$ are defined similarly. The matrix $A$ is composed of blocks $A_{ij}$, which has a block sparsity pattern, \emph{i.e.}, $A \in \mathbb{R}^{n \times n}(\mathcal{E}_p,0)$ with a partition $\{n_1, \ldots, n_N\}$ corresponding to the dimension of each subsystem's state. The matrices $B,M$ are of the forms $B = \text{diag}(B_1,\ldots,B_N)$ and $M = \text{diag}(M_1,\ldots,M_N)$. Our goal is to design a decentralized static state feedback
\begin{equation} \label{E:DeController}
    u_i(t) = -K_{ii} x_i(t), i = 1, \ldots, N
\end{equation}
such that the $\mathcal{H}_2$ norm of the transfer function $T_{zd}$ from disturbance $d$ to the desired  performance output $z$ is minimized. In~\eqref{E:DeController}, the global $K$ has a decentralized structure $\mathcal{K}$ as
\begin{equation*} 
    K \in \mathcal{K} := \{K \in \mathbb{R}^{m \times n} | K_{ij} = 0\; \text{if} \; i \neq j\},
\end{equation*}
where $m=\sum_{i=1}^N m_i, n = \sum_{i=1}^N n_i$, and each entry $K_{ij}$ is a block of dimension ${m_i \times n_j}$.

The design objective is
\begin{equation} \label{E:OptDeControl}
    \begin{aligned}
            \min_{K} \quad & \|T_{zd}\|^2 \\
            \text{s.t.} \quad & (A-BK) \; \text{is Hurwitz}, \\
            &  K \in \mathcal{K},
    \end{aligned}
\end{equation}
where $\|\cdot\|$ is the $\mathcal{H}_2$ norm of a transfer function. In this paper, the performance output $ z$ is chosen as
$$
    z = \begin{bmatrix} Q^{\frac{1}{2}} \\0 \end{bmatrix} x + \begin{bmatrix} 0\\ R^{\frac{1}{2}} \end{bmatrix} u,
$$
where $Q := \text{diag}(Q_1,\ldots, Q_N)$ and $R := \text{diag}(R_1,\ldots, R_N) $ denote the state and control performance weights, respectively, and diagonal block $Q_i, R_i$ correspond to the subsystem $i$. Adopting the same terminology in~\cite{siljak2011decentralized, geromel1994decentralized}, we refer to~\eqref{E:OptDeControl} as the \emph{optimal decentralized control} problem.

The constraint $\mathcal{K}$ does not allow any equivalent convex reformulation of the optimal decentralized problem~\eqref{E:OptDeControl} in general. Hence, problem~\eqref{E:OptDeControl} is challenging to solve exactly. Previous work either imposed special structures on system dynamics~\cite{rotkowitz2006characterization,tanaka2011bounded, shah2013optimal}, used certain relaxation/restriction techniques~\cite{geromel1994decentralized, fazelnia2017convex}, or applied non-convex optimization directly~\cite{lin2011augmented} to address this problem.

\subsection{Convex restriction via block-diagonal Lyapunov functions} \label{subsection:restriction}

It is well-known that the $\mathcal{H}_2$ norm of a stable linear system can be calculated using a linear matrix inequality~\cite{boyd1994linear}.
\begin{lemma}[\!\!\cite{boyd1994linear}]\label{lemma:H2norm}
    Consider a stable linear system $\dot{x}(t) = Ax(t)+Md(t), z(t) = Cx(t)$. The $\mathcal{H}_2$ norm of the transfer function from $d$ to $z$ can be computed by
    $$
        \|T_{zd}\|^2 = \inf_{X \succ 0}\{\textbf{Tr}\big(CXC^{\tr}\big) \mid AX+XA^{\tr}+MM^{\tr} \preceq 0\},
    $$
    where $\textbf{Tr}(\cdot)$ denotes the trace of a symmetric matrix.
\end{lemma}

According to Lemma~\ref{lemma:H2norm}, the optimal decentralized control problem~\eqref{E:OptDeControl} can be equivalently reformulated as
\begin{equation}\label{E:LMIHnorm1}
  \begin{aligned}
    \min_{X,K} \quad & {\textbf{Tr}}\big((Q+K^{\tr}RK)X\big) \\
    \text{s.t.} \quad & (A-BK)X+X(A-BK)^{\tr} + MM^{\tr} \preceq 0, \\
    & X \succ 0, K \in \mathcal{K}.
  \end{aligned}
\end{equation}
The first inequality in \eqref{E:LMIHnorm1} does not depend linearly on $X$ and $K$. A standard change of variables $ Z = KX $ leads to
\begin{equation*}
  \begin{aligned}
    \min_{X,Z} \quad & \textbf{Tr}(QX)+\textbf{Tr}\left(RZX^{-1}Z^{\tr}\right) \\
    \text{s.t.} \quad & (AX-BZ)+(AX-BZ)^{\tr} + MM^{\tr} \preceq 0, \\
    & X \succ 0, ZX^{-1} \in \mathcal{K}.
  \end{aligned}
\end{equation*}

To handle the nonlinear constraint $ZX^{-1} \in \mathcal{K}$, a classical parameterization idea~\cite{geromel1994decentralized} is to assume a \emph{block-diagonal} $X=\text{diag}(X_1,\ldots,X_N)$ with block size compatible to the subsystem's dimensions, which leads to 
$
    ZX^{-1} \in \mathcal{K} \Leftrightarrow Z \in \mathcal{K}.
$
 Considering the block-diagonal structures of $Q, R$, 
we have
$$
\begin{aligned}
    \textbf{Tr}(QX) &= \sum_{i=1}^N \textbf{Tr}\left(Q_iX_i\right), \\
    \textbf{Tr}\left(RZX^{-1}Z^{\tr}\right) &= \sum_{i=1}^N\textbf{Tr}\left(R_iZ_iX_i^{-1}Z_i^{\tr}\right).
\end{aligned}
$$
By introducing $Y_i \succeq Z_iX_i^{-1}Z_i^{\tr}$ and using the Schur complement~\cite{boyd1994linear}, a convex restriction to~\eqref{E:OptDeControl} is derived:
\begin{subequations}\label{E:LMIHnorm3}
  \begin{align}
    \min_{X_i,Y_i,Z_i} \;\; & \sum_{i=1}^N \textbf{Tr}(Q_iX_i)+\textbf{Tr}(R_iY_i) \nonumber \\
    \text{s.t.} \;\; & (AX-BZ)+(AX-BZ)^{\tr} + MM^{\tr} \preceq 0, \label{E:LMIHnorm3_eq1}\\
    & \begin{bmatrix} Y_i & Z_i \\ Z_i^{\tr} & X_i \end{bmatrix} \succeq 0, X_i \succ 0, i = 1, \ldots, N.
  \end{align}
\end{subequations}

Problem~\eqref{E:LMIHnorm3} is convex and ready to be solved using general conic solvers, and the decentralized controller is recovered as $K_{ii} = Z_iX_i^{-1}, i = 1, \ldots, N$.
In this paper, we make the following assumption.
\begin{assumption}
    Problem~\eqref{E:LMIHnorm3} is feasible, or equivalently system~\eqref{E:GlobalDynamics} is strongly decentralized stabilizable (see Definition~\ref{definition:strongly}).
\end{assumption}

\begin{remark} \label{Remark:Block-diagonalAssumption}
     The block-diagonal strategy was formally discussed in early 1990s~\cite{geromel1994decentralized}, which was later implicitly or explicitly used in the field of decentralized stabilization~\cite{siljak2011decentralized,Zheng2017Scalable,ahmadi2018distributed}. This strategy requires the closed-loop system to admit a block-diagonal Lyapunov function $V(x) = x^TPx = \sum_{i=1}^Nx_i^TP_ix_i$, where $P_i = X_i^{-1}, i \in \mathcal{V}$.  Problem~\eqref{E:LMIHnorm3} is a convex restriction of the original decentralized control problem~\eqref{E:OptDeControl}, and allows computing an upper bound of the optimal cost.
     However, quantifying the gap between the solution to~\eqref{E:LMIHnorm3} and the optimal solution to~\eqref{E:OptDeControl} is a challenging open problem, which is beyond the scope of this work. Indeed, problem~\eqref{E:LMIHnorm3} might be infeasible even for the cases in which problem~\eqref{E:OptDeControl} is feasible.
\end{remark}

\subsection{Problem statement}

   To connect the block-diagonal strategy with past work on decentralized control, we first present three classical definitions:
\begin{definition}[Stabilization]
    System~\eqref{E:GlobalDynamics} is called \emph{stabilizable}, if there exists a centralized controller $ u = -Kx $ such that the closed-loop system $\dot{x} = (A - BK)x$ is asymptotically stable.
\end{definition}

\begin{definition}[Decentralized stabilization~\cite{wang1973stabilization}]
    System~\eqref{E:GlobalDynamics} is called \emph{decentralized stabilizable}, if there exists a decentralized controller $ u_i = -K_{ii}x_i, i \in \mathcal{V}$ such that the closed-loop system $\dot{x} = (A - BK)x$ is asymptotically stable.
\end{definition}

\begin{definition}[Strongly decentralized stabilization~\cite{geromel1994decentralized}] \label{definition:strongly}
    System~\eqref{E:GlobalDynamics} is called \emph{strongly decentralized stabilizable} if there exists a decentralized $ u_i = -K_{ii}x_i, i \in \mathcal{V}$ such that the closed-loop system $\dot{x} = (A - BK)x$ admits a block-diagonal Lyapunov function $V(x) = \sum_{i=1}^Nx_i^TP_ix_i$.
\end{definition}

Then, we define three classes of complex systems:
\begin{equation*}
    \begin{aligned}
        \Sigma_0 &= \{(A,B) \mid \text{~\eqref{E:GlobalDynamics} is stabilizable}\}, \\
        \Sigma_1 &= \{(A,B) \mid \text{~\eqref{E:GlobalDynamics} is decentralized stabilizable}\}, \\
        \Sigma_2 &= \{(A,B) \mid \text{~\eqref{E:GlobalDynamics} is strongly decentralized stabilizable}\}. \\
    \end{aligned}
\end{equation*}
It is easy to see
$
    \Sigma_2 \subseteq \Sigma_1 \subseteq \Sigma_0.
$
In fact, the inclusion relationship is strict (see counterexamples in Appendix A), 
\begin{equation} \label{eq:Inclusion}
    \Sigma_2 \subset \Sigma_1 \subset \Sigma_0.
\end{equation}
The sets $\Sigma_0$ and $\Sigma_1$ can be algebraically characterized by {centralized fixed modes} and {decentralized fixed modes}~\cite{wang1973stabilization,alavian2014stabilizing}.  The class $\Sigma_2$ is useful for synthesizing decentralized controllers {as discussed in Section~\ref{subsection:restriction}}, but has been less studied before.
This motivates the first objective of our paper.

    \begin{problem}[Explicit characterizations]
        Derive sufficient conditions to characterize $\Sigma_2$.
    \end{problem}

   Solving~\eqref{E:LMIHnorm3} directly requires the global model knowledge, implicitly assuming the existence of a central entity to collect the complete model data. This performs a \emph{centralized design} of decentralized controllers. We note that the problem of distributed design using limited model information has received increasing attention~\cite{langbort2010distributed,farokhi2013optimal,ahmadi2018distributed, wang2017separable}. In our paper, we partition the subsystems into clusters to solve~\eqref{E:LMIHnorm3} in a distributed fashion. The second objective is as follows.

   \begin{problem}[Distributed computation]
       Given system~\eqref{E:GlobalDynamics} in $\Sigma_2$ with a plant graph $\mathcal{G}_p(\mathcal{V},\mathcal{E}_p)$, we aim 1) to partition the subsystems into {$t$ clusters $\mathcal{C}_1, \ldots, \mathcal{C}_t$}, where $\mathcal{V} = \cup_{i=1}^t\mathcal{C}_i$, and 2) to design a distributed algorithm to solve~\eqref{E:LMIHnorm3}, where the model data of subsystem $i$ is only shared within the clusters that contain it.
    \end{problem}

    We show that a chordal decomposition technique can be naturally used for {the partition $\mathcal{C}_1, \ldots, \mathcal{C}_t$}, and that the number of clusters depends on the sparsity of $\mathcal{G}_p(\mathcal{V},\mathcal{E}_p)$. We address Problem 1 in Section~\ref{Section:SufficientCondition} and Problem 2 in Sections~\ref{Section:Decomposition} and~\ref{Section:ADMM}.

\section{Sufficient Conditions on Strongly Decentralized Stabilization}\label{Section:SufficientCondition}

In this section, we discuss two classes of systems in $\Sigma_2$:
1) fully actuated systems, and 2) weakly coupled systems.

\subsection{Fully actuated systems}

\begin{definition}[Fully actuated systems]
    System~\eqref{E:GlobalDynamics} is called fully actuated, if each input matrix $B_{i}$ has full row rank, $i \in \mathcal{V}$.
\end{definition}

\begin{proposition} \label{proposition:fullyActuated}
    If system~\eqref{E:GlobalDynamics} is fully actuated, then we have $(A,B) \in \Sigma_2$.
\end{proposition}

\begin{IEEEproof}
    Consider the singular value decomposition of the input matrix $B_{i}$,
   \begin{equation*}
        B_{i} = U_i \begin{bmatrix} \Gamma_i & 0 \end{bmatrix}V_i^{\tr},
   \end{equation*}
   where $0$ is a zero block of  appropriate size, and $\Gamma_i  \in \mathbb{R}^{n_i \times n_i}$ is invertible since $B_{i}$ has full row rank. We then consider a decentralized feedback controller
   \begin{equation} \label{E:FullyDecen}
        K_{ii} = V_i \begin{bmatrix} \Gamma_i^{-1} \\ 0 \end{bmatrix}U_i^{\tr}(A_{ii} + \alpha_iI_{n_i}), i \in \mathcal{V},
   \end{equation}
   where $\alpha_i \in \mathbb{R}$. This choice leads to
   $$
        A_{ii} - B_{i}K_{ii} = -\alpha_i I_{n_i}, i \in \mathcal{V}.
   $$
    Using the decentralized controller~\eqref{E:FullyDecen}, the closed-loop system matrix becomes
    \begin{equation} \label{E:FullActuatedClosedLoop}
        \begin{aligned}
       A - BK
        =  \begin{bmatrix} -\alpha_1 I_{n_1} & A_{12} & \ldots & A_{1N} \\
                                A_{21} & -\alpha_2 I_{n_2} & \ldots & A_{2N} \\
                                \vdots & \vdots & \ddots & \vdots \\
                                A_{N1} & A_{N2} & \ldots & -\alpha_N I_{n_N} \\
                 \end{bmatrix}.
        \end{aligned}
    \end{equation}
    By choosing an appropriate $\alpha_i > 0$, we can always make $A - BK$ diagonally dominant with negative diagonal elements. Therefore,  $A - BK$ is diagonally stable, \emph{i.e.}, there exists a diagonal Lyapunov function to certify the stability of~\eqref{E:FullActuatedClosedLoop}. Therefore, we have $(A,B) \in \Sigma_2$.
\end{IEEEproof}

In essence, a fully actuated system is able to actuate each individual state directly. ~If each subsystem is of dimension one, \emph{i.e.}, $n_i = 1$, then the condition in Proposition~\ref{proposition:fullyActuated} means that the system pair of $(A_i,B_{i})$ is controllable. For general subsystems, the condition that $B_{i}$ has full row rank is stronger than the controllability of $(A_i,B_{i})$. Note that fully actuated systems have been used in some work on distributed design~\cite{langbort2010distributed, farokhi2013optimal}, where it required the input matrix $B$ to be invertible. Here, we show that a fully actuated system is indeed strongly decentralized stabilizable and suitable for the later development of the distributed algorithm.

\subsection{Weakly coupled systems}

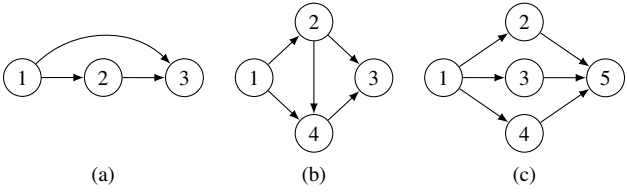
\begin{figure}[t]
    \centering
    \footnotesize
    \begin{tikzpicture}
	  \matrix (m) [matrix of nodes,
	  		       row sep = 0.8em,	
	  		       column sep = 2em,	
  			       nodes={circle, draw=black}] at (-2.8,0)
  		{ &  & \\ 1 & 2 & 3 \\ & &\\};
		\draw[-latex] (m-2-1) -- (m-2-2);
		\draw[-latex] (m-2-2) -- (m-2-3);
        \draw[-latex] (m-2-1)  to [out=45,in=135] (m-2-3);
		\node at (-2.8,-1.3) {(a)};
		\matrix (m2) [matrix of nodes,
	  		       row sep = 0.8em,	
	  		       column sep = 2em,	
  			       nodes={circle, draw=black}] at (2.8,0)
        { & 2 & \\ 1 & 3 & 5 \\& 4 &\\};
		\draw[-latex] (m2-2-1) -- (m2-2-2);
		\draw[-latex] (m2-2-1) -- (m2-1-2);
		\draw[-latex] (m2-2-1) -- (m2-3-2);
		\draw[-latex] (m2-1-2) -- (m2-2-3);
        \draw[-latex] (m2-3-2) -- (m2-2-3);
        \draw[-latex] (m2-2-2) -- (m2-2-3);
		\node at (0,-1.3) {(b)};

		\matrix (m3) [matrix of nodes,
	  		       row sep = 0.8em,	
	  		       column sep = 1.em,	
  			       nodes={circle, draw=black}] at (0,0)
        { & 2 & \\ 1 &  & 3 \\& 4 &\\};
		\draw[-latex] (m3-2-1) -- (m3-1-2);
		\draw[-latex] (m3-1-2) -- (m3-2-3);
		\draw[-latex] (m3-2-1) -- (m3-3-2);
		\draw[-latex] (m3-3-2) -- (m3-2-3);
        \draw[-latex] (m3-1-2) -- (m3-3-2);
		\node at (2.8,-1.3) {(c)};
	\end{tikzpicture}
    \caption{Examples of directed acyclic graphs}
    \label{F:uniGraph}
\end{figure}

Here, we discuss two types of weakly coupled systems: topologically weakly coupled systems and dynamically weakly coupled systems. A directed graph $\mathcal{G}$ is called \emph{acyclic} if there exist no directed cycles in $\mathcal{G}$. Fig.~\ref{F:uniGraph} shows some examples. A complex system with an acyclic $\mathcal{G}_p$ means that the dynamical influence among subsystems is unidirectional.

\begin{definition}[Topologically weakly coupled system]
    System~\eqref{E:GlobalDynamics} is called weakly coupled in terms of topological connections, if the plant graph $\mathcal{G}_p$ is acyclic.
\end{definition}

\begin{proposition} \label{proposition:TopoWeak}
    For the class of topologically weakly coupled systems, we have
    $$\Sigma_1 = \Sigma_2 = \{(A,B) \mid (A_{ii}, B_{i}) \text{~is stabilizable}, i \in \mathcal{V}\}.$$
\end{proposition}

\begin{IEEEproof}
    This result is a simple consequence of~\cite{sootla2017block,carlson1992block}.  
    If $\mathcal{G}_p$ is acyclic, then there exists an ordering of the nodes such that for every edge $(v_1,v_2)$, node $v_1$ precedes node $v_2$ in the ordering. For this ordering, the resulting system matrix $A$ is block lower triangular. Thus, without loss of generality, for a topologically weakly coupled system~\eqref{E:GlobalDynamics}, the closed-loop system with a decentralized controller remains block lower triangular.    It is known that a block triangular matrix is stable if and only if it is block-diagonally stable~\cite{sootla2017block,carlson1992block}, \emph{i.e.}, there exists a block-diagonal Lyapunov function to certify the stability of the closed-loop system. Therefore, for the class of topologically weakly coupled systems, we have
    $
        (A,B) \in \Sigma_1 \Leftrightarrow (A,B) \in \Sigma_2.
    $
    Meanwhile, considering the block triangular structure, the overall closed-loop system is stable if and only if each isolated closed-loop subsystem $A_{ii} - B_{i}K_{ii}$ is stable, $i \in \mathcal{V}$. This completes the proof.
\end{IEEEproof}

We note that the class of topologically weakly coupled systems is also known as  \emph{hierarchical systems}; see~\cite[Chapter 10]{lunze1992feedback}. Hierarchical systems have useful properties, \emph{e.g.}, the equivalence between stability and block-diagonal stability~\cite{sootla2017block,carlson1992block}. Proposition~\ref{proposition:TopoWeak} further shows that for this type of systems, decentralized stabilization is equivalent to strongly decentralized stabilization ($\Sigma_1 = \Sigma_2$).

Next, we consider dynamically weakly coupled systems. If each pair $(A_{ii},B_{i})$ is stabilizable, then there exists a local feedback $K_{ii}$ such that $A_{ii} - B_{i}K_{ii}$ is stable. Consequently, given any $Q_i \succ 0$, there exists a $P_i \succ 0$, such that
$$
    (A_{ii} - B_{i}K_{ii})^{\tr}P_i + P_i(A_{ii} - B_{i}K_{ii}) + Q_i \prec 0.
$$
In some cases, \emph{e.g.}, the singular values of $A_{i j}$ are small (the strength of interactions is low), there may still exist a solution $P_i \succ 0$ for the following inequality
\begin{equation} \label{E:WeakCoupling0}
    \begin{aligned}
    (A_{ii} - B_{i}K_{ii})^{\tr}P_i &+ P_i(A_{ii} - B_{i}K_{ii}) \\
        &+ P_i\bigg(\sum_{j\in \mathbb{N}_i} A_{ij}A_{ij}^{\tr}\bigg)P_i + Q_i \prec 0.
    \end{aligned}
\end{equation}
In~\eqref{E:WeakCoupling0}, recall that $\mathbb{N}_i$ denotes  the set of neighbouring nodes of node $i$. This observation leads to a concept of dynamically weakly coupled systems.
\begin{definition}[Dynamically weakly coupled systems] \label{Defi:DynaWeak}
    System~\eqref{E:GlobalDynamics} is weakly coupled in terms of dynamical interactions, if there exists a local feedback $K_{ii}$ such that the following inequality holds
    \begin{equation} \label{E:WeakCoupling}
        \begin{aligned}
        (A_{ii} &- B_{i}K_{ii})^{\tr}P_i + P_i(A_{ii} - B_{i}K_{ii})  \\
        &+ P_i\left(\sum_{j\in \mathbb{N}_i} A_{ij}W_{ij}^{-1}A_{ij}^{\tr}\right)P_i + \sum_{j\in \hat{\mathbb{N}}_i} W_{ji} \prec 0,
        \end{aligned}
    \end{equation}
    for some $W_{ij} \succ 0, j \in {\mathbb{N}}_i, P_i \succ 0, i \in \mathcal{V}$, where 
    $\hat{\mathbb{N}}_i$ denotes the set of nodes coming out of node $i$ in $\mathcal{G}_p$.
\end{definition}

Definition~\ref{Defi:DynaWeak} is more general than condition~\eqref{E:WeakCoupling0}, since inequality~\eqref{E:WeakCoupling} is reduced to~\eqref{E:WeakCoupling0} when setting $W_{ij} = I_{n_j}, j \in \mathbb{N}_i$, and $Q_i = \sigma_i I_{n_i}$, where $\sigma_i$ denotes the number of nodes in $\hat{\mathbb{N}}_i$.

\begin{proposition} \label{Prop:dynamicallyWeaklyCoupled}
    For a dynamically weakly coupled system~\eqref{E:GlobalDynamics}, \emph{i.e.},~\eqref{E:WeakCoupling} holds, we have $(A,B) \in \Sigma_2$.
\end{proposition}

The proof utilizes the following Lemma.

\begin{lemma} \label{lemma:LMI}
    Given two matrices $X,Y$ of appropriate dimensions, we have 
    \begin{equation} \label{eq:lemma1}
        X^{\tr}WX + Y^{\tr}W^{-1}Y \succeq X^{\tr}Y + Y^{\tr}X
     \end{equation}
     for any $W \succ 0$ of appropriate dimension.
\end{lemma}
\begin{IEEEproof}
    Observe that
    \begin{equation*}
        \begin{aligned}
            &X^{\tr}WX + Y^{\tr}W^{-1}Y - (X^{\tr}Y+Y^{\tr}X) \\
            = \;\; &(WX-Y)^{\tr} W^{-1}(WX - Y) \succeq 0.
        \end{aligned}
    \end{equation*}
    This means~\eqref{eq:lemma1} holds.
\end{IEEEproof}

  \textit{Proof of Proposition~\ref{Prop:dynamicallyWeaklyCoupled}:}  Consider a decentralized controller
    $
        K = \text{diag}(K_{11}, \ldots, K_{NN}).
    $
     Upon defining $\hat{A}_{ii} = A_{ii} - B_{i}K_{ii}$ and ignoring the disturbance, the closed-loop dynamics for each subsystem become
    \begin{equation}\label{eq:closedSub}
         \dot{x}_i(t) = \hat{A}_{ii}x_i(t) + \sum_{j\in \mathbb{N}_i}A_{ij}x_j(t), \quad \forall \; i \in \mathcal{V}.
    \end{equation}
    Next, we consider a block-diagonal Lyapunov function
    $
        V(x) = \sum_{i=1}^N x_i^{\tr}(t)P_ix_i(t).
    $
    The derivative of $V(x)$ along the closed-loop trajectory~\eqref{eq:closedSub} is
    \begin{equation} \label{eq:DeriLyapunov}
        \begin{aligned}
            \dot{V}(x) &= \sum_{i=1}^N \big(\dot{x}_i^{\tr}P_ix_i + x_i^{\tr}P_i\dot{x}_i\big) \\
                       &=  \sum_{i=1}^N \bigg(x_i^{\tr}\left(\hat{A}_{ii}^{\tr}P_i+P_i\hat{A}_{ii}\right)x_i + \\
                      & \qquad \underbrace{\big(\sum_{j\in \mathbb{N}_i}A_{ij}x_j\big)^{\tr}P_ix_i + x_i^{\tr}P_i\big(\sum_{j\in \mathbb{N}_i}A_{ij}x_j\big)}_{\text{coupling term}}\bigg).
        \end{aligned}
    \end{equation}
    For the coupling term in~\eqref{eq:DeriLyapunov}, according to {Lemma~\ref{lemma:LMI}}, we have
    \begin{equation} \label{eq:CouplingTerm}
        \begin{aligned}
             &\big(\sum_{j\in \mathbb{N}_i}A_{ij}x_j\big)^{\tr}P_ix_i + x_i^{\tr}P_i\big(\sum_{j\in \mathbb{N}_i}A_{ij}x_j\big) \\
=\:& \sum_{j\in \mathbb{N}_i} \big(x_j^{\tr}A_{ij}^{\tr}P_ix_i + x_i^{\tr}P_iA_{ij}x_j\big) \\
            \leq \: &\sum_{j\in \mathbb{N}_i} \big(x_i^{\tr}P_iA_{ij}W_{ij}^{-1}A_{ij}^{\tr}P_ix_i + x_j^{\tr}W_{ij}x_j\big),
        \end{aligned}
    \end{equation}
    for any $W_{ij} \succ 0, j \in \mathbb{N}_i$.   Substituting~\eqref{eq:CouplingTerm} into \eqref{eq:DeriLyapunov}, we get
    \begin{equation*} 
    \small
        \begin{aligned}
            \dot{V}(x) \leq &\sum_{i=1}^N \bigg(x_i^{\tr}\big(\hat{A}_{ii}^{\tr}P_i+P_i\hat{A}_{ii}  \\
            & \qquad +    P_i(\sum_{j \in \mathbb{N}_i}A_{ij}W_{ij}^{-1}A_{ij}^{\tr})P_i\big)x_i
            + \sum_{j\in \mathbb{N}_i}x_j^TW_{ij}x_j \bigg)\\
            = & \sum_{i=1}^N x_i^{\tr}\bigg(\hat{A}_{ii}^TP_i+P_i\hat{A}_{ii} \\
            & \qquad +
             P_i\big(\sum_{j\in \mathbb{N}_i}A_{ij}W_{ij}^{-1}A_{ij}^{\tr}\big)P_i  + \sum_{j\in \hat{\mathbb{N}}_i} W_{ji}\bigg)x_i.
              \end{aligned}
    \end{equation*}
    If condition~\eqref{E:WeakCoupling} holds  for some $W_{ij} \succ 0, j \in \mathbb{N}_i$, $P_i \succ 0, i \in \mathcal{V}$,
    then, $\dot{V}(x)$ is negative definite. Thus, $V(s)$ is a block-diagonal Lyapunov function for the closed-loop system.
    \qed

Note that condition~\eqref{E:WeakCoupling} can be equivalently formulated into the following problem:
we aim to find static scaling matrices $W_{ij} \succ 0$ such that there exist $K_{ii}, i \in \mathcal{V}$ satisfying
\begin{equation} \label{Eq:dissipativity}
    \|\hat{W}_i (sI - A_{ii} + B_{i}K_{ii})^{-1}\hat{A}_iW_i\|_{\infty} < 1, \;\; i \in \mathcal{V},
\end{equation}
where
$
    \hat{A}_i = \begin{bmatrix} A_{ij_1} & A_{ij_2}& \ldots& A_{ij_s} \end{bmatrix},   W_i = \text{diag}\{W_{ij_1}^{-\frac{1}{2}},W_{ij_2}^{-\frac{1}{2}}, \ldots ,W_{ij_s}^{-\frac{1}{2}}\}, \hat{W}_i = \left(\sum_{j \in \hat{\mathbb{N}}_i} W_{ji}\right)^{\frac{1}{2}}
$ and $\mathbb{N}_i = \{j_1, \ldots, j_s\}$. It is clear that both~\eqref{E:WeakCoupling} and ~\eqref{Eq:dissipativity} are coupled between subsystems due to the scaling matrices $W_{ij} \succ 0$. If we \emph{a priori} fix the weights $W_{ij}$,  then the constraints in~\eqref{E:WeakCoupling} and~\eqref{Eq:dissipativity} are decoupled. This leads to a set of localized conditions to certify the dynamically weakly coupled condition~\eqref{E:WeakCoupling}. The sufficient conditions for block-diagonal stability based on scaled diagonal
dominance in~\cite{sootla2017block} may be good choices for choosing the weights $W_{ij}$.

\section{Chordal decomposition in optimal decentralized control} \label{Section:Decomposition}

In this section, by assuming that an undirected version of the plant graph $\mathcal{G}_p(\mathcal{V},\mathcal{E}_p)$ is chordal, we derive a decomposed version of problem~\eqref{E:LMIHnorm3}, leading to multiple local subproblems. The chordal structure provides a way to define local computing agents or clusters of subsystems. This facilitates us to develop a distributed algorithm to solve~\eqref{E:LMIHnorm3} in Section~\ref{Section:ADMM}.

\subsection{Chordal graphs and sparse matrices}

For completeness, we first review some graph-theoretic notion, and refer the interested reader to~\cite{blair1993introduction, vandenberghe2014chordal} for details. Graph $\mathcal{G}(\mathcal{V},\mathcal{E})$ is called undirected if $(i,j) \in \mathcal{E} \Leftrightarrow (j,i) \in \mathcal{E}$.  A \emph{clique} $\mathcal{C}$ is a subset of nodes in $\mathcal{V}$ where any pair of distinct nodes has an edge, \emph{i.e.}, $(i,j) \in \mathcal{E}, \forall i,j \in \mathcal{C}, i \neq j$. If a clique $\mathcal{C}$ is not included in any other clique, then it is called a \emph{maximal clique}. A \emph{cycle} of length $k$ is a sequence of nodes $\{v_1, v_2, \ldots, v_k\} \subseteq \mathcal{V}$ with $(v_k, v_{1}) \in \mathcal{E}$ and $(v_i, v_{i+1}) \in \mathcal{E}, \forall i = 1, \ldots, k-1$. A \emph{chord} in a cycle $\{v_1, v_2, \ldots, v_k\}$ is an edge $(v_i,v_j)$ that joins two non-adjacent nodes in the cycle.

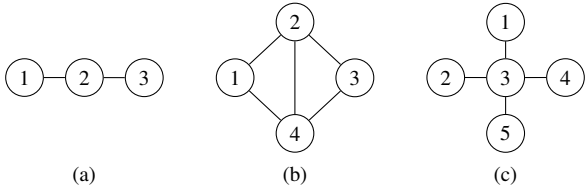
\begin{figure}[t]
    \centering
    \footnotesize
    \begin{tikzpicture}
	  \matrix (m) [matrix of nodes,
	  		       row sep = 0.8em,	
	  		       column sep = 1em,	
  			       nodes={circle, draw=black}] at (-2.8,0)
  		{ &  & \\ 1 & 2 & 3 \\ & &\\};
		\draw (m-2-1) -- (m-2-2);
		\draw (m-2-2) -- (m-2-3);
		\node at (-2.8,-1.3) {(a)};
		\matrix (m2) [matrix of nodes,
	  		       row sep = 0.8em,	
	  		       column sep = 1em,	
  			       nodes={circle, draw=black}] at (2.8,0)
        { & 1 & \\ 2 & 3 & 4 \\& 5 &\\};
		\draw (m2-1-2) -- (m2-2-2);
		\draw (m2-2-1) -- (m2-2-2);
		\draw (m2-2-2) -- (m2-2-3);
		\draw (m2-2-2) -- (m2-3-2);
		\node at (0,-1.3) {(b)};

		\matrix (m3) [matrix of nodes,
	  		       row sep = 0.8em,	
	  		       column sep = 1.em,	
  			       nodes={circle, draw=black}] at (0,0)
        { & 2 & \\ 1 &  & 3 \\& 4 &\\};
		\draw (m3-1-2) -- (m3-2-1);
		\draw (m3-1-2) -- (m3-2-3);
		\draw (m3-2-1) -- (m3-3-2);
		\draw (m3-2-3) -- (m3-3-2);
        \draw (m3-1-2) -- (m3-3-2);
		\node at (2.8,-1.3) {(c)};
	\end{tikzpicture}
    \caption{Examples of chordal graphs: (a) a line graph; (b) a triangulated graph; (c) a star graph.}
    \label{F:ChordalGraph}
    \vspace{-3mm}
\end{figure}

An undirected graph $\mathcal{G}$ is called \emph{chordal} if every cycle of length at least four has one chord~\cite{blair1993introduction}. Note that the set of maximal cliques is unique in a chordal graph, and the graph decomposition based on the maximal cliques is unique accordingly~\cite{vandenberghe2014chordal}. Fig.~\ref{F:ChordalGraph} illustrates some examples, and there are two maximal cliques, $\mathcal{C}_1 = \{1, 2\}$ and $\mathcal{C}_2 = \{2, 3\}$ for the chordal graph shown in Fig.~\ref{F:ChordalGraph}(a). We highlight that maximal cliques can serve as computing agents and the overlapping elements, \emph{e.g.}, node 2 in Fig.~\ref{F:ChordalGraph}(a), will play a role of coordination among maximal cliques. This feature enables preserving model data privacy (see Remarks~\ref{R:Privacy} and~\ref{R:Privacy&Cliques}). 

Given a sequence of integers $\{n_1, \ldots, n_N\}$, and an undirected graph $\mathcal{G}(\mathcal{V},\mathcal{E})$, we define the space of symmetric block matrices with a particular sparsity pattern as
\begin{equation*} 
    \mathbb{S}^n(\mathcal{E},0) := \{X \in \mathbb{S}^{n} | X_{ij} = X_{ji}^{\tr}= 0\; \text{if} \; (j,i) \notin \hat{\mathcal{E}} \},
\end{equation*}
where $n = \sum_{i=1}^N n_i$, $X_{ij}$ is a block of dimension $n_i \times n_j$ and $\hat{\mathcal{E}} = \mathcal{E} \cup \{(i,i), \forall i \in \mathcal{V}\}$. The cone of sparse block PSD matrices is defined as
\begin{equation*}
    \mathbb{S}^n_+(\mathcal{E},0) := \{X \in \mathbb{S}^{n}(\mathcal{E},0) | X \succeq 0 \}.
\end{equation*}

Given a partition $\{n_1,\ldots, n_N\}$ and a maximal clique $\mathcal{C}_k$ of $\mathcal{G}$, we define a block index matrix $E_{\mathcal{C}_k} \in \mathbb{R}^{|\mathcal{C}_k| \times n}$ with $|\mathcal{C}_k| = \sum_{j \in \mathcal{C}_k}n_j$ and $n = \sum_{i=1}^Nn_i$ as
$$
    (E_{\mathcal{C}_k})_{ij} := \begin{cases} I_{n_i}, \qquad \;\, \text{if } \mathcal{C}_k(i) = j, \\ 0_{n_i \times n_j}, \quad \text{otherwise}, \end{cases}
$$
where $\mathcal{C}_k(i)$ denotes the $i$-th node in $\mathcal{C}_k$, sorted in the natural ordering, $I_{n_i}$ denotes an identity matrix of size $n_i \times n_i$, and $0_{n_i \times n_j}$ denotes a zero matrix of size $n_i \times n_j$. Note that $X_k = E_{\mathcal{C}_k}XE_{\mathcal{C}_k}^{\tr} \in \mathbb{S}^{|\mathcal{C}_k|}$ extracts a principal submatrix according to clique $\mathcal{C}_k$, and the operation $E_{\mathcal{C}_k}^{\tr}X_kE_{\mathcal{C}_k}$ {inflates} a $\vert\mathcal{C}_k\vert \times \vert \mathcal{C}_k\vert $ matrix into a sparse $n\times n$ matrix. Then, we have the following result.

\begin{lemma}[\!\!\!\cite{griewank1984existence,agler1988positive,ZKSP2017scalable}]\label{T:ChordalDecompositionTheorem}
     Let $\mathcal{G}(\mathcal{V},\mathcal{E})$ be a chordal graph with {maximal cliques $\{\mathcal{C}_1,\mathcal{C}_2, \ldots, \mathcal{C}_t\}$}. Given a partition $\{n_1,n_2,\ldots,n_N\}$, we have $X\in\mathbb{S}^n_{+}(\mathcal{E},0)$ if and only if there exist matrices $X_k \in \mathbb{S}^{\vert \mathcal{C}_k \vert}_+, k=1,\,\ldots,\,t$, such that
    $$X = \sum_{k=1}^{t} E_{\mathcal{C}_k}^{\tr} X_k E_{\mathcal{C}_k}.$$
\end{lemma}

\begin{example}
    Consider the following positive semidefinite matrix with a trivial partition $\{1,1,1\}$
    \begin{equation*} 
        \begin{bmatrix} 2 & 1 & 0 \\ 1 & 1 & 1 \\ 0 & 1 & 2 \end{bmatrix} \succeq 0,
    \end{equation*}
    which has a chordal sparsity pattern corresponding to Fig~\ref{F:ChordalGraph}(a) with maximal cliques $\mathcal{C}_1 = \{1,2\}$ and $\mathcal{C}_2 = \{2,3\}$. Then, Lemma~\ref{T:ChordalDecompositionTheorem} guarantees the following decomposition 
    $$
         \begin{bmatrix} 2 & 1 & 0 \\ 1 & 1 & 1 \\ 0 & 1 & 2 \end{bmatrix} =  E_{\mathcal{C}_1}^\tr\underbrace{\begin{bmatrix} 2 & 1  \\ 1 & 0.5  \end{bmatrix}}_{\succeq 0}E_{\mathcal{C}_1}  + E_{\mathcal{C}_2}^\tr\underbrace{\begin{bmatrix} 0.5 & 1 \\  1 & 2 \end{bmatrix}}_{\succeq 0}E_{\mathcal{C}_2}.
    $$
    where $$
    E_{\mathcal{C}_1} = \begin{bmatrix} 1&0&0 \\0& 1& 0  \end{bmatrix}, E_{\mathcal{C}_2} = \begin{bmatrix} 0&1&0 \\0& 0& 1  \end{bmatrix}.
    $$
    Indeed, for any PSD matrix with partition $\{2,1,1\}$ corresponding to Fig.~\ref{F:ChordalGraph}(a), Lemma~\ref{T:ChordalDecompositionTheorem} guarantees a block-wise decomposition as follows ($*$ denotes a real number)
     \begin{equation*} 
        \begin{aligned}
       & \underbrace{\left[
	\begin{array}{cc|c|cc}
	* & * & * & 0   \\
	* & * & * & 0   \\\hline
	* & * & * & * \\\hline
    0 & 0 & * & * \\
	\end{array}
	\right]}_{\succeq 0} 
 =   \underbrace{\left[
	\begin{array}{cc|c|c}
	* & * & * & 0   \\
	* & * & * & 0   \\\hline
	* & * & * & 0 \\\hline
    0 & 0 & 0 & 0 \\
	\end{array}\right]}_{\succeq 0} +    \underbrace{\left[
	\begin{array}{cc|c|c}
	0 & 0 & 0 & 0   \\
	0 & 0 & 0 & 0   \\\hline
	0 & 0 & * & * \\\hline
    0 & 0 & * & * \\
	\end{array}\right]}_{\succeq 0}.
\end{aligned}
    \end{equation*}
\end{example}

\subsection{Chordal decomposition of problem~\eqref{E:LMIHnorm3}} \label{subsection:ChordalDecomposition}

In~\eqref{E:LMIHnorm3}, the variables $X_i,Y_i,Z_i$ are coupled by the inequality~\eqref{E:LMIHnorm3_eq1} only, while the rest of the constraints and the objective function are naturally separable due to the separable performance weights $Q,R$. Meanwhile, thanks to the block-diagonal assumption on $X$, the coupled linear matrix inequality has a structured sparsity pattern characterized by an undirected version of $\mathcal{G}_p(\mathcal{V},\mathcal{E}_p)$. Precisely, we define an undirected graph
$
    \mathcal{G}_u(\mathcal{V}, \mathcal{E}_u)
$
with $\mathcal{E}_u = \mathcal{E}_p \cup \mathcal{E}_p^{\tr}$, where $\mathcal{E}_p^{\tr}$ denotes the edge set of the transpose graph of $\mathcal{G}_p$, \emph{i.e.}, the graph associated to the transpose of the adjacency matrix of $\mathcal{G}_p$.
%
\begin{assumption}
    Graph $\mathcal{G}_u$ is chordal with maximal cliques $\mathcal{C}_1, \ldots, \mathcal{C}_p$.
\end{assumption}

\begin{figure}[t]
    \centering
    \footnotesize
    \begin{tikzpicture}
	  \matrix (m) [matrix of nodes,
	  		       row sep = 0.8em,	
	  		       column sep = 2em,	
  			       nodes={circle, draw=black}] at (-2,0)
  		{ 1 & 2 & 3 \\};
		\draw[-latex] (m-1-1) -- (m-1-2);
		\draw[-latex] (m-1-2) -- (m-1-3);
		\node at (-2,-0.7) {(a)};

		\matrix (m3) [matrix of nodes,
	  		       row sep = 0.8em,	
	  		       column sep = 2.em,	
  			       nodes={circle, draw=black}] at (2,0)
      { 1 & 2 & 3 \\};	
      \draw[-latex] (m3-1-2) -- (m3-1-1);
		\draw[-latex] (m3-1-3) -- (m3-1-2);
		\node at (2,-0.7) {(b)};
	\end{tikzpicture}
    \caption{Constructing the transpose graph: (a) directed plant graph $\mathcal{G}_p(\mathcal{V},\mathcal{E}_p)$, (b) the transpose graph of $\mathcal{G}_p$. }
    \label{F:TransposeGraph}
\end{figure}
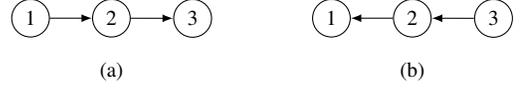

\begin{remark}
    The undirected graph $\mathcal{G}_u$ will be used in the development of distributed computation using ADMM. For example, consider an interconnected system with a directed line graph in Fig.~\ref{F:TransposeGraph}(a). Its transpose graph is shown in Fig.~\ref{F:TransposeGraph}(b), and the resulting undirected graph $\mathcal{G}_u$ is the same as that in Fig.~\ref{F:ChordalGraph}(a). If $\mathcal{G}_u$ is not chordal, we can add suitable edges to $\mathcal{E}_u$ to obtain a chordal graph~\cite{vandenberghe2014chordal}. In this case, sharing model data with directed neighbours in $\mathcal{G}_p$ is not sufficient for the proposed distributed solution. Still, privacy of model data is maintained within each maximal clique in $\mathcal{G}_u$.  For simplicity, we assume that $\mathcal{G}_u$ is chordal. As shown in Fig.~\ref{F:ChordalGraph}, some graphs, such as chains, trees and stars, are already chordal. 
\end{remark}

Considering the inherent structure of system~\eqref{E:SubDynamics}, it is straightforward to see that
 $   (AX-BZ)+(AX-BZ)^{\tr} + MM^{\tr} \in \mathbb{S}^n(\mathcal{E}_u,0).$
To ease the exposition, we define
\begin{equation*}
    F(X,Z) := - (AX-BZ) - (AX-BZ)^{\tr} - MM^{\tr}.
\end{equation*}
According to Lemma~\ref{T:ChordalDecompositionTheorem}, $F(X,Z) \succeq 0 $ is equivalent to the condition that there exist $J_k \in \mathbb{S}^{|\mathcal{C}_k|}_{+}, k = 1, \ldots, t$, such that
\begin{equation} \label{E:H2DeConstraint}
    F(X,Z) = \sum_{k=1}^t E_{\mathcal{C}_k}^{\tr} J_k E_{\mathcal{C}_k}.
\end{equation}

Therefore,~\eqref{E:LMIHnorm3} can be equivalently decomposed into
\begin{equation}\label{E:H2Decomposition}
  \begin{aligned}
    \min_{X_i,Y_i,Z_i,J_k} \quad & \sum_{i=1}^N \textbf{Tr}(Q_iX_i)+\textbf{Tr}(R_iY_i) \\
    \text{s.t.~~~} \quad & \sum_{k=1}^t E_{\mathcal{C}_k}^{\tr} J_k E_{\mathcal{C}_k} = F(X,Z), \\
    & \begin{bmatrix} Y_i & Z_i \\ Z_i^{\tr} & X_i \end{bmatrix} \succeq 0, X_i \succ 0, i = 1, \ldots, N, \\
    & J_k \succeq 0, k = 1, \ldots, t. \\
  \end{aligned}
\end{equation}
One notable feature of~\eqref{E:H2Decomposition} is that the global constraint~\eqref{E:LMIHnorm3_eq1} is replaced by a set of small coupled constraints~\eqref{E:H2DeConstraint}. In other words,~\eqref{E:H2Decomposition} has partially coupled constraints, which can be solved in a distributed way by introducing consensus variables.

The cliques $\mathcal{C}_1, \ldots, \mathcal{C}_t$ give a partition of subsystems, and will serve as local computing agents. If there is no overlap among the cliques $\mathcal{C}_1, \ldots, \mathcal{C}_t$ (\emph{i.e.}, the system~\eqref{E:GlobalDynamics} is composed by dynamically disjoint components), then~\eqref{E:H2Decomposition} is trivially decomposed into $t$ decoupled subproblems of decentralized optimal control, which can be solved by cliques $\mathcal{C}_1, \ldots, \mathcal{C}_t$ independently. In the case where different cliques share some common nodes with each other, we can introduce appropriate auxiliary variables to achieve a distributed solution using ADMM.

\section{Distributed Design based on ADMM} \label{Section:ADMM}

{To formulate our distributed approach to solve the decomposed problem~\eqref{E:H2Decomposition} (equivalent to problem~\eqref{E:LMIHnorm3}), we briefly review the basic setup of ADMM; see~\cite{boyd2011distributed} for a comprehensive review.
ADMM is a first-order method that solves a convex optimization problem of the form
\begin{equation}
\label{E:ADMM_base_prob}
    \begin{aligned}
        \min_{x,y} \quad & f(x)+g(y) \\
        \text{s.t.} \quad & Ex + Fy = c,
    \end{aligned}
\end{equation}
where $x \in \mathbb{R}^{n_x}, y \in \mathbb{R}^{n_y}$ are decision variables, $f$ and $g$ are convex functions, and $ E \in \mathbb{R}^{n_c\times n_x}, F \in \mathbb{R}^{n_c\times n_y}$ and $c \in \mathbb{R}^{n_c}$ are problem data. Given a penalty parameter $\rho>0$, the scaled ADMM algorithm solves~\eqref{E:ADMM_base_prob} using the following iterations
\begin{equation*}
    \begin{aligned}
        x^{h+1} & = \text{arg} \min_{x} \, f(x) + \frac{\rho}{2}\|Ex + Fy^h - c + \lambda^h\|^2,
        \\
        y^{h+1} & = \text{arg} \min_{y} \, g(y) + \frac{\rho}{2}\|Ex^{h+1} + Fy - c + \lambda^h\|^2,
        \\
        \lambda^{h+1} &= \lambda^{h} +  E x^{h+1} + F y^{h+1} - c, 
    \end{aligned}
\end{equation*}
where $\lambda \in \mathbb{R}^{n_c}$ is a scaled dual variable, and $h$ denotes the iteration index. In many applications, splitting the minimization over $x$ and $y$ often leads to multiple subproblems, allowing distributed computation; see~\cite{boyd2011distributed} for detailed discussions.

\vspace{1mm}
\subsection{A simple example} \label{Section:SimpleExample}

To illustrate the approach, we first consider an interconnected system characterized by a chain of three nodes, as shown in Fig.~\ref{F:ChordalGraph}(a). In this case, the model data are $B = \text{diag}\{B_1,B_2,B_3\}, M = \text{diag}\{M_1,M_2,M_3\}$ and
$$
    A = \begin{bmatrix} A_{11} & A_{12} & 0 \\ A_{21} & A_{22} & A_{23} \\ 0 & A_{32}& A_{33} \end{bmatrix}.
$$
Note that the following illustration is directly suitable for systems with a directed graph. For example, if the plant graph is a directed line as in~Fig.~\ref{F:TransposeGraph}(a), then we have $A_{12} = 0, A_{23} = 0$ in matrix $A$ and need to construct the same undirected graph in Fig.~\ref{F:ChordalGraph}(a) for the distributed computation.

In this case, there are two cliques $\mathcal{C}_1 = \{1,2\}, \mathcal{C}_2 = \{2,3\}$, and $J_k, k = 1,2$ in~\eqref{E:H2DeConstraint} are in the following form
\begin{equation*} 
    \begin{aligned}
        J_1(X_1,X_2,Z_1,J_{22,1}) &:= -\begin{bmatrix} J_{11} & A_{12}X_2 + X_1A_{21}^{\tr} \\ * & J_{22,1}\end{bmatrix}, \\
        J_2(X_2,X_3,Z_3,J_{22,2}) &:= -\begin{bmatrix} J_{22,2} & A_{23}X_3 + X_2A_{32}^{\tr} \\ * & J_{33}\end{bmatrix},
    \end{aligned}
\end{equation*}
where $*$ denotes the corresponding symmetric part and
\begin{equation*} 
    \begin{aligned}
        J_{11} := A_{11}X_1 - B_1Z_1 + (A_{11}X_1 - B_1Z_1)^{\tr} + M_1M_1^{\tr}, \\
        J_{33} := A_{33}X_3 - B_3Z_3 + (A_{33}X_3 - B_3Z_3)^{\tr} + M_3M_3^{\tr}. 
    \end{aligned}
\end{equation*}
The coupling effect is imposed on the overlapping node 2:
\begin{equation*}
    \begin{aligned}
        J_{22,1} + J_{22,2} = &A_{22}X_2 - B_2Z_2  \\ &+ (A_{22}X_2 - B_2Z_2)^{\tr} + M_2M_2^{\tr}.
    \end{aligned}
\end{equation*}

For any coupling variables that appear in two cliques, we introduce auxiliary variables. For this case, we introduce auxiliary variables for node 2
\begin{equation} \label{E:H2Exauxiliary}
    \begin{aligned}
        J_{22,1} &= \hat{J}_{22,1}, J_{22,2} = \hat{J}_{22,2}, \\
        X_{2}   &= X_{2,1}, \quad X_{2}   = X_{2,2}.
    \end{aligned}
\end{equation}
Also, we split the variables according to the cliques and the overlapping node
\begin{equation*}
    \begin{aligned}
        \text{Node 2, } \quad \;\, y &:= \{X_{2},Y_2,Z_2,\hat{J}_{22,1},\hat{J}_{22,2}\}, \\
        \text{Clique $\mathcal{C}_1$, } \;\; x_{\mathcal{C}_1} &:= \{X_1,Y_1,Z_1,X_{2,1},J_{22,1}\}, \\
        \text{Clique $\mathcal{C}_2$, } \;\; x_{\mathcal{C}_2} &:=\{X_3,Y_3,Z_3,X_{2,2},J_{22,2}\}. \\
    \end{aligned}
\end{equation*}
The variable $y$ corresponds to the same $y$ in the canonical form~\eqref{E:ADMM_base_prob} and variables $x_{\mathcal{C}_1}, x_{\mathcal{C}_2}$ corresponds to $x$ in~\eqref{E:ADMM_base_prob}. This can be seen more directly in~\eqref{E:H2Example}.
Next, we show that~\eqref{E:H2Decomposition} can be rewritten into the standard ADMM form~\eqref{E:ADMM_base_prob} by defining indicator functions as
\begin{equation*}
    \begin{aligned}
    &\mathbb{I}_{\mathcal{S}_k}(x_{\mathcal{C}_k}) := \begin{cases} 0, \qquad x_{\mathcal{C}_k} \in \mathcal{S}_k, \\ +\infty, \quad \text{otherwise}, \end{cases} \\
    &\mathbb{I}_{\mathcal{L}}(y) := \begin{cases} 0, \qquad y_l \in \mathcal{L}, \\ +\infty, \quad \text{otherwise}, \end{cases}
    \end{aligned}
\end{equation*}
where sets ${\mathcal{S}}_1, {\mathcal{S}}_2$ are defined as
$$
\begin{aligned}
    &\begin{aligned}
        {\mathcal{S}}_1 :=  \bigg\{ x_{\mathcal{C}_1} \bigg\vert   J_1(X_1,X_{2,1},&Z_1,J_{22,1}) \succeq 0,  X_1 \succ 0,  \\
        &\begin{bmatrix} Y_1 & Z_1 \\ Z_1^{\tr} & X_1 \end{bmatrix} \succeq 0 \text{~are feasible}
        \bigg\},
    \end{aligned} \\
    &\begin{aligned}
        {\mathcal{S}}_2 :=  \bigg\{ x_{\mathcal{C}_2} \bigg\vert   J_2(X_{2,2},X_3,&Z_3,J_{22,2}) \succeq 0,  X_3 \succ 0,  \\
        &\begin{bmatrix} Y_3 & Z_3 \\ Z_3^{\tr} & X_3 \end{bmatrix} \succeq 0 \text{~are feasible}
        \bigg\},
    \end{aligned}
\end{aligned}
$$
and ${\mathcal{L}}$ is defined by
$$
    \begin{aligned}
        {\mathcal{L}} :=  \bigg\{ y \bigg\vert  \hat{J}_{22,1} &+ \hat{J}_{22,2} = A_{22}X_2 - B_2Z_2  \\
                          &+ (A_{22}X_2 - B_2Z_2)^{\tr} + M_2M_2^{\tr},   \\
        & \quad X_2 \succ 0, \begin{bmatrix} Y_2 & Z_2 \\ Z_2^{\tr} & X_2 \end{bmatrix} \succeq 0 \text{~are feasible}
        \bigg\}.
    \end{aligned} \\
$$

This allows us to rewrite~\eqref{E:H2Decomposition} as an optimization problem in the form of~\eqref{E:ADMM_base_prob}:
\begin{equation}\label{E:H2Example}
  \begin{aligned}
    \min_{x_{\mathcal{C}_k,y}} \quad & \sum_{k=1}^2 f_k(x_{\mathcal{C}_k}) + g(y) \\
    \text{s.t.~} \quad &  \text{\eqref{E:H2Exauxiliary} holds},
  \end{aligned}
\end{equation}
where $f_1(x_{\mathcal{C}_1}),f_2(x_{\mathcal{C}_2}$) based on each clique are defined as
\begin{subequations} \label{E:H2ExFunction}
    \begin{align}
       f_1(x_{\mathcal{C}_1})&:= \textbf{Tr}(Q_1X_1)+\textbf{Tr}(R_1Y_1) + \mathbb{I}_{{\mathcal{S}}_1}(x_{\mathcal{C}_1}), \label{E:H2ExFunction_s1}\\
       f_2(x_{\mathcal{C}_2})&:= \textbf{Tr}(Q_3X_3)+\textbf{Tr}(R_3Y_3) + \mathbb{I}_{{\mathcal{S}}_2}(x_{\mathcal{C}_2}), \label{E:H2ExFunction_s1}
    \end{align}
\end{subequations}
and $g(y)$ based on the overlapping node 2 is defined as
\begin{equation*}
    g(y):= \textbf{Tr}(Q_2X_2)+\textbf{Tr}(R_2Y_2) + \mathbb{I}_{{\mathcal{L}}}(y).
\end{equation*}

Upon denoting $\hat{x}_{\mathcal{C}_k}$ as the variables in $x_{\mathcal{C}_k}$ that appears in the consensus constraint~\eqref{E:H2Exauxiliary}, and $y_l({\mathcal{C}_k})$ as the corresponding local copies, \emph{e.g.},
$
    \hat{x}_{\mathcal{C}_1} = \{X_{2,1}, J_{22,1}\}, y({\mathcal{C}_1}) =  \{X_{2}, \hat{J}_{22,1}\},
$
the ADMM algorithm for~\eqref{E:H2Example} takes a distributed form:
\begin{enumerate}
    \item[] \emph{ADMM algorithm for the distributed design}
  \item $x$-update: for each clique $k$, solve the local problem:
  \begin{equation}\label{E:ExampleH2_S1}
  x_{\mathcal{C}_k}^{h+1} = \arg \min_{x_{\mathcal{C}_k}} f_k(x_{\mathcal{C}_k}) + \frac{\rho}{2}\|\hat{x}_{\mathcal{C}_k} - y^h({\mathcal{C}_k}) + \lambda_{\mathcal{C}_k}^h\|^2.
  \end{equation}
  \item $y$-update: solve the following problem to update local variables
  \begin{equation} \label{E:ExampleH2_S2}
   y^{h+1} = \arg \min_{y} g(y)+ \frac{\rho}{2}\sum_{k=1}^2\|\hat{x}_{\mathcal{C}_k}^{h+1} - y(\mathcal{C}_k) + \lambda_{\mathcal{C}_k}^h\|^2.
    \end{equation}
  \item $\lambda$-update: compute the dual variable
  \begin{equation} \label{E:ExampleH2_S3}
  \lambda_{\mathcal{C}_k}^{h+1} = \lambda_{\mathcal{C}_k}^{h} + \hat{x}_{\mathcal{C}_k}^{h+1} - y^{h+1}({\mathcal{C}_k}), k = 1,2.
   \end{equation}
\end{enumerate}

\begin{figure}[t]
    \centering
    \setlength{\belowcaptionskip}{0em}
    \small
    \begin{tikzpicture}
        \matrix (m1) [matrix of nodes,
	  		       row sep = 1.8em,	
	  		       column sep = 3em,	
  			       nodes={circle, draw=black},
                   column 2/.style={nodes={fill=blue!20}}] at (0,-2.)
        { 1 & 2 &  \\ & 2 & 3 \\};
		\draw (m1-1-1) -- (m1-1-2);
        \draw (m1-2-2) -- (m1-2-3);
        \draw[blue,thick,dashed,<->] (0,-1.7) -- (0,-2.2);
        \node at (-0.9,-1.7) {$\mathcal{C}_1$};
        \node at (0.9,-2.3) {$\mathcal{C}_2$};
	\end{tikzpicture}
    \caption{ Illustration of the ADMM algorithm for solving~\eqref{E:H2Example}: cliques $\mathcal{C}_1 = \{1, 2\}$ and $\mathcal{C}_2 = \{2, 3\}$ serve as two computing agents and the overlapping node 2 plays a role of a coordinator. 
    }
    \label{F:H2Example}
    \vspace{-3mm}
\end{figure}
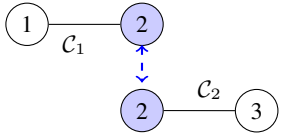

At each iteration $h$, subproblem~\eqref{E:ExampleH2_S1} only depends on each clique $\mathcal{C}_k$. Consequently, the cliques can serve as computing agents to solve subproblem~\eqref{E:ExampleH2_S1} to update the variable $x_{\mathcal{C}_k}^{h+1}$ in parallel. For example, clique $\mathcal{C}_1$ needs to solve the following convex problem
$$
    \begin{aligned}
             \min_{x_{\mathcal{C}_1}} \;\; & \textbf{Tr}(Q_1X_1)+\textbf{Tr}(R_1Y_1) + \frac{\rho}{2}\|\hat{x}_{\mathcal{C}_1} - y^h({\mathcal{C}_1}) + \lambda_{\mathcal{C}_1}^h\|^2 \\
             \text{s.t.~} \;\; & \begin{bmatrix} J_{11} & A_{12}X_{2,1} + X_1A_{21}^{\tr} \\ * & J_{22,1}\end{bmatrix} \preceq 0, \\
             & \begin{bmatrix} Y_1 & Z_1 \\ Z_1^{\tr} & X_1 \end{bmatrix} \succeq 0, X_1 \succ 0.
   \end{aligned}
$$
where the regularization term is
$$
    \begin{aligned}
        \|\hat{x}_{\mathcal{C}_1} - y^h({\mathcal{C}_1}) + \lambda_{\mathcal{C}_1}^h\|^2  = \|X_{2,1} &- X_2^h + \Lambda_{2,1}^h\|^2 \\
        &+ \|J_{22,1} - \hat{J}_{22,1}^h + \Lambda_{22,1}^h\|^2.
    \end{aligned}
$$
The subproblems~\eqref{E:ExampleH2_S2} and~\eqref{E:ExampleH2_S3} deal with the consensus variables $y^{h+1}$ and multipliers $\lambda_{\mathcal{C}_k}^{h+1}, k = 1,2$, which can be computed by node 2. Fig.~\ref{F:H2Example} illustrates the distributed nature of this algorithm.

\begin{remark}[Privacy of model data] \label{R:Privacy}
    At each iteration, the coordinator (\emph{i.e.}, node 2) only requires model data of itself $A_{22}, B_2, M_2$ and the local copies $X_{2,k}^{h+1}, J_{22,k}^{h+1}$ from cliques $\mathcal{C}_k, k=1,2$ to update $y^{h+1}, \lambda_{\mathcal{C}_k}^{h+1}$ by solving~\eqref{E:ExampleH2_S2} and~\eqref{E:ExampleH2_S3}. Therefore, the proposed ADMM algorithm for solving~\eqref{E:H2Decomposition} has a distributed nature (see Fig.~\ref{F:H2Example} for illustration): cliques $\mathcal{C}_1$ and $\mathcal{C}_2$ can solve~\eqref{E:ExampleH2_S1} based on \emph{the model data within each clique} in parallel, and node 2 plays a role of coordination by updating the auxiliary variables $y^{h+1}, \lambda_{\mathcal{C}_k}^{h+1}$. Consequently, the model data of node $1$ (\emph{i.e.}, $A_{11},B_1,M_1,A_{12},A_{21}$) are accessible only to clique $\mathcal{C}_1$ only, while clique $\mathcal{C}_2$ holds the model data of node $3$ (\emph{i.e.}, $A_{33},B_3,M_3,A_{32},A_{23}$), \emph{exclusively}.
\end{remark}

\begin{remark}[Privacy and maximal cliques]\label{R:Privacy&Cliques}
    In our ADMM algorithm, the privacy of model data are maintained within each maximal clique of $\mathcal{G}_u$. Therefore, the level of privacy depends on the sparsity of $\mathcal{G}_u$. For highly interconnected systems with only one maximal clique, the decomposition~\eqref{E:H2DeConstraint} brings no benefit for privacy, and a global model is still required. In practice, if the plant graph $\mathcal{G}_p$ is a chain or star graph (see Fig.~\ref{F:ChordalGraph} for illustration), then each maximal clique is of size two only, meaning that each subsystem need to share its model data with its direct neighbors only, and the model data privacy can be therefore maintained to a large extent.
\end{remark}

\begin{remark}[Convergence of the ADMM algorithm]
    The general ADMM algorithm is guaranteed to converge for convex problems under very mild conditions~\cite[Section 3.2]{boyd2011distributed}. In our case, under the feasibility assumption of~\eqref{E:LMIHnorm3}, the proposed ADMM algorithm~\eqref{E:ExampleH2_S1}-\eqref{E:ExampleH2_S3} is guaranteed to find a solution asymptotically. 
    {In the examples considered in this work}, ADMM typically found a solution with moderate accuracy (in the sense of standard stopping criteria~\cite[Section 3.3]{boyd2011distributed}) within a few hundred iterations (see Section~\ref{Section:Results}). Note that adjusting the penalty parameter $\rho$ dynamically may further improve the practical convergence of the ADMM algorithm~\cite[Section 3.4.1]{boyd2011distributed}. In our simulations, we used a fixed choice of $\rho = 5$, since it led to a satisfactory convergence for our instances.
\end{remark}

\subsection{The general case}

The idea above can be extended to solve~\eqref{E:H2Decomposition} with a general chordal graph pattern, and the general problem~\eqref{E:H2Decomposition} shares great similarities with the simple example in Section~\ref{Section:SimpleExample}. First, we define a set $\mathcal{N}_0 := \{i \in \mathcal{V} \mid \exists q, k = 1,\ldots,p, \text{such that~} i \in \mathcal{C}_q \cap \mathcal{C}_k\}$ that contains the overlapping nodes, and a set $\mathcal{E}_0 := \{(i,j) \in \mathcal{E}_u \mid \exists q, k = 1,\ldots,t, \text{such that~} (i,j) \in (\mathcal{C}_q\times\mathcal{C}_q) \cap (\mathcal{C}_k\times\mathcal{C}_k)\}$ that contains the overlapping edges. For the example in~Fig.~\ref{F:H2Example}, we have $\mathcal{N}_0 = \{2\} $ and $\mathcal{E}_0 = \emptyset$. Also, we define $\mathcal{N}_i := \{k \mid i \in \mathcal{C}_k, k = 1,\ldots, p\}$ that denotes the cliques containing node $i$, and $\mathcal{E}_{ij} := \{ k \mid (i,j) \in \mathcal{C}_k \times \mathcal{C}_k, k = 1,\ldots, t\}.$ that denotes the cliques containing edge $(i,j)$.

In fact, the elements in $\mathcal{N}_0 $ and $\mathcal{E}_0$ make the constraint~\eqref{E:H2DeConstraint} coupled among different maximal cliques. Similar to~\eqref{E:H2Exauxiliary}, for each node $i \in \mathcal{N}_0$, we introduce  local
consensus constraints
\begin{equation} \label{E:NodeConsensus}
    X_i = X_{i,k}, \hat{J}_{ii,k}=J_{ii,k}, \forall k \in \mathcal{N}_i.
\end{equation}
For each overlapping edge $(i,j) \in \mathcal{E}_0$, we introduce local consensus constraints
\begin{equation} \label{E:EdgeConsensus}
    X_{ij} = X_{i,k}, \hat{J}_{ij,k}= J_{ij,k}, \forall (i,j) \in \mathcal{E}_{ij}.
\end{equation}
Then, variable $x_{\mathcal{C}_k}$ for each maximal clique $k=1,\ldots,t$ includes
\begin{itemize}
  \item $X_i,Y_i,Z_i, i \in \mathcal{C}_k\setminus \mathcal{N}_0$ that belongs to clique $\mathcal{C}_k$ exclusively;
  \item $X_{i,k}, J_{ii,k}, i \in \mathcal{C}_k \cap \mathcal{N}_0$ that corresponds to overlapping nodes in $\mathcal{C}_k$ ;
  \item $J_{ij,k}, (i,j) \in (\mathcal{C}_k\times \mathcal{C}_k) \cap \mathcal{E}_0$ that corresponds to overlapping edges in $\mathcal{C}_k$;
\end{itemize}
We also collect the local copies $X_i,Y_i, Z_i,\hat{J}_{ii,k}, i \in \mathcal{N}_0$ and $\hat{J}_{ij,k}, X_{ij,k}, (i,j) \in\mathcal{E}_0$ as the consensus variable $y$.

Then,~\eqref{E:H2Decomposition} can be written into the canonical ADMM form:
\begin{equation}\label{E:H2General}
  \begin{aligned}
    \min_{x_{\mathcal{C}_k,y}} \quad & \sum_{k=1}^t f_k(x_{\mathcal{C}_k}) + g(y) \\
    \text{s.t.} \quad & \text{\eqref{E:NodeConsensus} and \eqref{E:EdgeConsensus} hold},
  \end{aligned}
\end{equation}
where the clique function $f_k(x_{\mathcal{C}_k})$ is defined as
\begin{equation} \label{E:CliqueFunction}
    f_k(x_{\mathcal{C}_k}) := \sum_{i \in \mathcal{C}_k\setminus \mathcal{N}_0} \textbf{Tr}(Q_iX_i)+\textbf{Tr}(R_iY_i) + \mathbb{I}_{{\mathcal{S}}_k}(x_{\mathcal{C}_k}),
\end{equation}
and $g(y)$ is defined as
\begin{equation} \label{E:OverlappingFunction}
    g(y):= \sum_{i \in \mathcal{N}_0} \textbf{Tr}(Q_iX_i)+\textbf{Tr}(R_iY_i) + \mathbb{I}_{\mathcal{L}}(y).
\end{equation}
In~\eqref{E:CliqueFunction}, set $\mathcal{S}_k$ is defined as
$$
\begin{aligned}
        {\mathcal{S}}_k :=  \bigg\{ x_{\mathcal{C}_k} \bigg\vert  & J_k(x_{\mathcal{C}_k}) \succeq 0,  X_i \succ 0,  \\
        &\begin{bmatrix} Y_i & Z_i \\ Z_i^{\tr} & X_i \end{bmatrix} \succeq 0, i \in \mathcal{C}_k\setminus \mathcal{N}_0 \text{~are feasible}
        \bigg\},
\end{aligned}
$$
and in~\eqref{E:OverlappingFunction}, set $\mathcal{L}$ is defined as \vspace{-2mm}
$$
{\small
    \begin{aligned}
        {\mathcal{L}} :=   \bigg\{ y \bigg\vert  &\sum_{k\in \mathcal{N}_i} \hat{J}_{ii,k} = A_{ii}X_i - B_iZ_i + (A_{ii}X_i - B_iZ_i)^{\tr} \\
        &\quad+ M_iM_i^{\tr},  X_i \succ 0, \begin{bmatrix} Y_i & Z_i \\ Z_i^{\tr} & X_i \end{bmatrix} \succeq 0, i \in \mathcal{N}_0, \\
        &  \sum_{k\in \mathcal{E}_{ij}} \hat{J}_{ij,k} = A_{ij}X_{ij} + X_{ji}A_{ji}^{\tr},
        (i,j) \in \mathcal{E}_0 \text{~are feasible}
        \bigg\}.
    \end{aligned}
}
$$%
\vspace{-3 mm}

By applying the ADMM to~\eqref{E:H2General}, we obtain iterations that are identical to~\eqref{E:ExampleH2_S1}-\eqref{E:ExampleH2_S3}. Note that the  set 
$\mathcal{L}$ can be equivalently rewritten as a product of sets defined by $X_i,Y_i, Z_i,\hat{J}_{ii,k}, i \in \mathcal{N}_0$ and $\hat{J}_{ij,k}, X_{ij}, (i,j) \in\mathcal{E}_0$. For each $i \in \mathcal{N}_0$, the set for $X_i,Y_i, Z_i,\hat{J}_{ii,k}$ is defined as
$$
    \begin{aligned}
        \mathcal{L}_i :&=\bigg\{ (X_i,Y_i, Z_i,\hat{J}_{ii,k}) \bigg\vert  \sum_{k\in \mathcal{N}_i} \hat{J}_{ii,k} = A_{ii}X_i - B_iZ_i +\\
        &  (A_{ii}X_i - B_iZ_i)^{\tr} + M_iM_i^{\tr}, X_i \succ 0, \begin{bmatrix} Y_i & Z_i \\ Z_i^{\tr} & X_i \end{bmatrix} \succeq 0 \bigg\}. 
    \end{aligned}
$$
This means that $y$-update~\eqref{E:ExampleH2_S2} can be distributed among the overlapping nodes $\mathcal{N}_0$ and overlapping edges $\mathcal{E}_0$. Therefore, similar to the example in Section~\ref{Section:SimpleExample}, variables $x_{\mathcal{C}_k}^{h}$ can be updated on each clique in parallel, and the overlapping elements in $\mathcal{N}_0 $ and $\mathcal{E}_0$ can update $y_{\mathcal{C}_k}^{h}, \lambda_{\mathcal{C}_k}^{h}$ individually until convergence.

Here, as stated in Remark~\ref{R:Privacy}, we emphasize that the main interest of our algorithm is the ability of distributing the computation to cliques and overlapping elements, thus preserving the privacy of model data in the problem.

\section{Numerical Cases} \label{Section:Results}

This section demonstrates the effectiveness of the proposed distributed design method\footnote{\scriptsize Code is available via \url{https://github.com/zhengy09/distributed_design_methods}.}. For the examples, we ran the ADMM algorithm with termination tolerance $10^{-3}$ and the number of iterations was limited to 500.
In our simulations, SeDuMi~\cite{sturm1999using} and YALMIP~\cite{lofberg2004yalmip} were used to solve the subproblems within each clique and overlapping elements.

\subsection{First-order systems with acyclic directed graphs}

\begin{figure}
    \centering
    \small
    \begin{tikzpicture}
        \matrix (m1) [matrix of nodes,
	  		       row sep = 0.9em,	
	  		       column sep = 1em,	
  			       nodes={circle, draw=black},
                   column 2/.style={nodes={fill=blue!20}}] at (-1.3,0)
        { & 2 \\ 1 & \\& 4 \\};
		\draw (m1-1-2) -- (m1-2-1);
		\draw (m1-1-2) -- (m1-3-2);
		\draw (m1-2-1) -- (m1-3-2);
        \draw[blue,thick,dashed,<->] (-0.5,0.85) -- (0.5,0.85);
        \draw[blue,thick,dashed,<->] (-0.5,-0.85) -- (0.5,-0.85);
        \node at (-1.2,0) {$\mathcal{C}_1$};

        \matrix (m2) [matrix of nodes,
	  		       row sep = .9em,	
	  		       column sep = 1.em,	
  			       nodes={circle, draw=black},
                   column 1/.style={nodes={fill=blue!20}}] at (1.3,0)
        { 2 & \\   & 3 \\ 4 & \\};
		\draw (m2-1-1) -- (m2-2-2);
		\draw (m2-1-1) -- (m2-3-1);
		\draw (m2-2-2) -- (m2-3-1);
        \node at (1.2,0) {$\mathcal{C}_2$};
	\end{tikzpicture}
    \caption{Illustration of the ADMM algorithm for solving~\eqref{E:LMIHnorm3} corresponding to the example~\eqref{E:ExampleEx1}: the cliques $\mathcal{C}_1 = \{1,2,4\}$ and $\mathcal{C}_2 = \{2,3,4\}$ serve as two computing agents and the overlapping nodes play a role of coordinators by updating the axillary variables.}
    \label{F:ExampleCoordinate}
    \vspace{-4mm}
\end{figure}
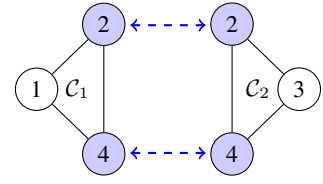

We first consider a network of four unstable coupled first-order subsystems, where $\mathcal{G}_p$ is the directed acyclic graph shown in Fig.~\ref{F:uniGraph}(b). In the experiment, the global dynamics are
\begin{equation} \label{E:ExampleEx1}
    \dot{x}(t) = \begin{bmatrix} 1 & 0 & 0 & 0 \\
                              1 & 2 & 0 & 0 \\
                              0 & 2 & 3 & 4 \\
                              1 & 2 & 0 & 4   \end{bmatrix} x(t) + u(t) + d(t).
\end{equation}
This system is both fully actuated and topologically weakly coupled according to Section~\ref{Section:SufficientCondition}. 
We chose $Q_i = 1$ and $R_i = 1, i \in \mathcal{V}$ in our simulation. When the global dynamics are available, solving~\eqref{E:LMIHnorm3} directly returned a decentralized controller $ K_{11} = 7.34; K_{22} = 11.38; K_{33} = 6.16, K_{44} = 13.48$ with an $\mathcal{H}_2$ performance of 5.36.

Instead, when the privacy of model data is concerned, the proposed ADMM algorithm can solve~\eqref{E:LMIHnorm3} in a distributed fashion. As shown in Fig.~\ref{F:ExampleCoordinate}, for clique 1, only the model data of nodes 1, 2, 4 are required, while clique 2 only needs the model data of nodes 2, 3, 4, and the overlapping nodes 2 and 4 play a role of coordinations in the algorithm. In this way, the model of node 1 can be kept private within clique 1 and the model of node 3 is known within clique 2 exclusively. For this instance, after 54 iterations, the ADMM algorithm returned the decentralized controller $K_{11} = 7.35; K_{22} = 11.41; K_{33} = 6.16, K_{44} = 13.49$ with an $\mathcal{H}_2$ performance of 5.37. The convergence plot of our algorithm for this instance is given in Fig.~\ref{fig:residual}.
\begin{figure}[t]%
\setlength{\abovecaptionskip}{3pt}
\centering
\subfigure[]{%
\label{fig:first}%
\includegraphics[height=1.4in]{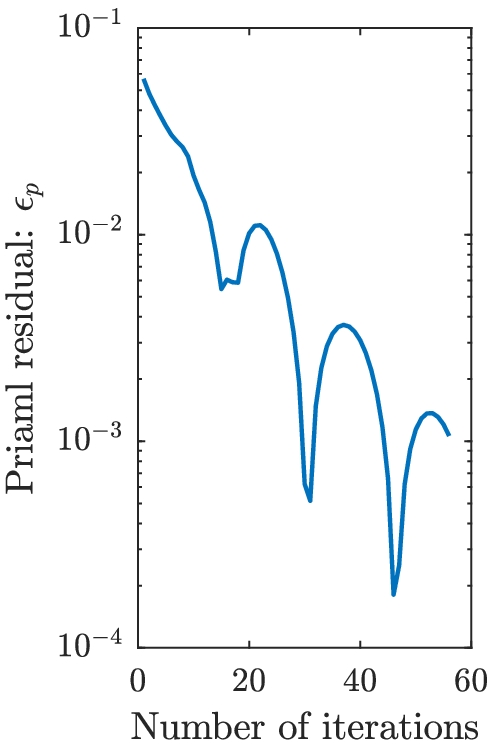}}%
\qquad
\subfigure[]{%
\includegraphics[height=1.4in]{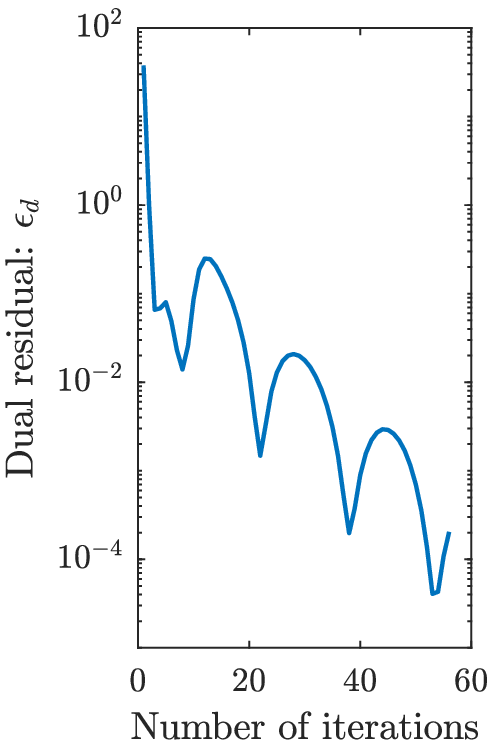}}%
\caption{Primal and dual residuals versus iteration number for the directed graph (42): (a) primal residual, (b) dual residual.}
\label{fig:residual}
\end{figure}

\subsection{A chain of unstable second-order coupled systems} \label{Section:comparsion}

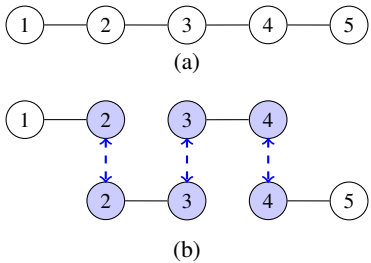
\begin{figure}[t]
    \centering
    \setlength{\abovecaptionskip}{0em}
    \setlength{\belowcaptionskip}{0em}
    \footnotesize
    \begin{tikzpicture}
        \matrix (m) [matrix of nodes,
	  		       row sep = 1em,	
	  		       column sep = 2em,	
  			       nodes={circle, draw=black},
                   ] at (0,0)
        { 1 & 2 & 3 & 4 & 5\\};
		\draw (m-1-1) -- (m-1-2);
		\draw (m-1-2) -- (m-1-3);
        \draw (m-1-3) -- (m-1-4);
        \draw (m-1-4) -- (m-1-5);
        \node at (0,-0.5) {\small (a)};

        \matrix (m1) [matrix of nodes,
	  		       row sep = 2em,	
	  		       column sep = 2em,	
  			       nodes={circle, draw=black},
                   column 2/.style={nodes={fill=blue!20}},
                   column 3/.style={nodes={fill=blue!20}},
                   column 4/.style={nodes={fill=blue!20}}] at (0,-1.8)
        { 1 & 2 & 3 & 4 &  \\ & 2 & 3 & 4& 5 \\ };
		\draw (m1-1-1) -- (m1-1-2);
        \draw (m1-2-2) -- (m1-2-3);
         \draw (m1-1-3) -- (m1-1-4);
          \draw (m1-2-4) -- (m1-2-5);
        \node at (0,-3) {\small (b)};	
        \draw[blue,thick,dashed,<->] (-1.08,-1.5) -- (-1.08,-2.1);
        \draw[blue,thick,dashed,<->] (0,-1.5) -- (0,-2.1);
        \draw[blue,thick,dashed,<->] (1.08,-1.5) -- (1.08,-2.1);

	\end{tikzpicture}
    \caption{(a) A chain of five nodes;
    (b) Four maximal cliques $\mathcal{C}_i = \{i,i+1\}, i = 1, 2, 3, 4$, which serve as four computing agents relying only on the model data within each clique; the overlapping nodes 2, 3, 4 play a role of coordinators.}
    \label{F:ExChain}
    \vspace{-4mm}
\end{figure}


Here, we use a chain of five nodes (see Fig.~\ref{F:ExChain}) to provide a comparison between the proposed ADMM algorithm and the following three approaches: 
\begin{enumerate}
  \item \emph{A sequential approach}~\cite{Zheng2017Scalable}, which exploits the properties of clique trees in chordal graphs;
  \item \emph{Localized LQR design}~\cite[Chapter 7.3]{lunze1992feedback}, which computes a local LQR controller for each subsystem independently by ignoring the coupling terms $A_{ij}$;
  \item \emph{Truncated LQR design}, which computes a centralized LQR controller using the global model data and only keeps the diagonal blocks for decentralized feedback.
\end{enumerate}

\begin{figure}[t]
    \centering
        \setlength{\abovecaptionskip}{0.2em}
    \includegraphics[width=0.5\columnwidth]{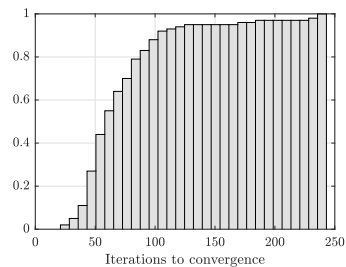}
    \caption{Cumulative plot of the fraction of 100 random trails of~\eqref{E:SimEx1} that required a given number of iterations to converge.} 
    \label{F:ChainIterations}
    \vspace{-1mm}
\end{figure}

It is assumed that each node is an unstable second order system coupled with its neighbouring nodes,
\begin{equation} \label{E:SimEx1}
{
    \dot{x}_i = \begin{bmatrix} 1 & 1 \\ 1 & 2 \end{bmatrix} x_i + \sum_{j \in \mathbb{N}_i}A_{ij}x_j + \begin{bmatrix}0 \\1\end{bmatrix} (u_i + d_i),}
\end{equation}
where the entries of coupling term $A_{ij}$ were generated randomly from $-0.5$ to $0.5$ to ensure that the numerical examples are strongly decentralized stabilizable. There are four maximal cliques $\mathcal{C}_i = \{i,i+1\}, i = 1, 2, 3, 4$. 
The model data can be kept private within each clique, and the overlapping nodes (\emph{i.e.}, $2, 3, 4$) coordinate the consensus variables among maximal cliques. %
In the simulation, the state and control weights were  $Q_i = I_2$ and $R_i = 1$ for each subsystem.

We generated 100 random instances of this interconnected system~\eqref{E:SimEx1}. The performance comparison between the four methods is listed in Table~\ref{Table:Comparision}. The proposed ADMM algorithm was able to return stabilizing decentralized controllers for all 100 tests, while the sequential approach, localized LQR and truncated LQR design only succeeded for 72\,\%, 54\,\%, 62\,\% of the tests, respectively. This is expected since the proposed ADMM algorithm only requires the system being strongly decentralized stabilizable. 
The sequential approach requires an additional equal-splitting assumption among maximal cliques (see~\cite[Section VI.B]{Zheng2017Scalable}), and the localized LQR and truncated LQR design has no guarantees of success in general. 
Also, the average $\mathcal{H}_2$ performance for the common succeeded instances by the ADMM algorithm is the best. Finally, Fig.~\ref{F:ChainIterations} shows the cumulative plot of convergence performance of our algorithm, where 90\,\% of the tests required less than 150 iterations.

\begin{table}[t]
    \setlength{\abovecaptionskip}{0.2em}
    \caption{Comparison of the proposed ADMM algorithm, sequential approach~\cite{Zheng2017Scalable}, localized LQR~\cite[Chapter 7.3]{lunze1992feedback} and truncated LQR design for the system~\eqref{E:SimEx1}.}	
    \centering
    \setlength{\abovecaptionskip}{3pt}
    \setlength{\belowcaptionskip}{0em}
    \label{Table:Comparision}
    \begin{tabular*}{\linewidth}{@{\extracolsep{\fill}}ccccc}
         \toprule
          &  ADMM &  Sequential &  {Localized LQR} & {Truncated LQR}\\
         \midrule
         pct.$^{\ddagger}$ & 100\,\% & 72\,\% & 54\,\% & 62\,\%\\
         $\mathcal{H}_2$$^{\dagger}$ & 6.06 & 6.36 & 6.50 & 6.49\\
         \bottomrule
       \end{tabular*}
       \scriptsize
   \newline
   \raggedright
   $^\ddagger$: Successful percentage of returning a stabilizing decentralized controller. \\
   $\,^\dagger$: Average $\mathcal{H}_2$ performance of 
   based on common successful instances.\\
   \vspace{-4mm}
\end{table}

\balance
\section{Conclusion} \label{Section:Conclusion}

We introduced a distributed design method for decentralized control that relies on local model information only. Our main strategy is consistent with the recent general idea of exploiting sparsity in systems theory via chordal decomposition~\cite{mason2014chordal, andersen2014robust, ZKSP2017scalable, Zheng2017Scalable}. In this paper, we further demonstrated the potential of chordal decomposition in distributed design of decentralized controllers, by combining this approach with the ADMM algorithm. Similar to~\cite{Zheng2017Scalable,ahmadi2018distributed}, 
our method relies on a block-diagonal Lyapunov function, which may bring some conservatism in general.
{Currently, we are 
studying convex restrictions that are less restrictive than the block-diagonal assumption, while still allowing distributed computation.}

\appendix

\section{Counterexamples}

This appendix shows~\eqref{eq:Inclusion} using counterexamples. Consider the following system with two scalar subsystems
\begin{equation} \label{E:Example}
    \begin{bmatrix} \dot{x}_1 \\ \dot{x}_2 \end{bmatrix} = \begin{bmatrix} 1 & 2\\ a_1 & a_2 \end{bmatrix}  \begin{bmatrix} {x}_1 \\ {x}_2 \end{bmatrix} +  \begin{bmatrix} 0 & 0 \\ 0 & 1 \end{bmatrix}  \begin{bmatrix} u_1 \\ u_2 \end{bmatrix},
\end{equation}
where the first scalar subsystem is not affected by the first control input, \emph{i.e.}, $B_{1} = 0$ in~\eqref{E:SubDynamics}.
Since~\eqref{E:Example} is controllable $\forall a_1 \in \mathbb{R}, a_2 \in \mathbb{R}$, then we know
$$
    \left( \begin{bmatrix} 1 & 2\\ a_1 & a_2 \end{bmatrix}, \begin{bmatrix} 0 & 0 \\ 0 & 1 \end{bmatrix}\right) \in \Sigma_0, \forall a_1 \in \mathbb{R}, a_2 \in \mathbb{R}.
$$

Next, consider a decentralized controller for~\eqref{E:Example}
$
    u_1 = -k_1 x_1, u_2 = -k_2 x_2,
$
then the closed-loop system becomes
\begin{equation} \label{E:ExampleClosedloop}
    \begin{bmatrix} \dot{x}_1 \\ \dot{x}_2 \end{bmatrix} = \begin{bmatrix} 1 & 2\\ a_1 & a_2-k_2 \end{bmatrix}  \begin{bmatrix} {x}_1 \\ {x}_2 \end{bmatrix}.
\end{equation}
The stability of~\eqref{E:ExampleClosedloop} means that the real parts of its eigenvalues are negative. This requires 
$$
    \begin{cases}
        a_2 + 1 - k_2 < 0, \\
        a_2 - k_2 - 2a_1 > 0,
    \end{cases}
$$
which is equivalent to
$
    a_2 + 1 <  k_2 < a_2 - 2a_1.
$
This means
$$
    \left( \begin{bmatrix} 1 & 2\\ a_1 & a_2 \end{bmatrix}, \begin{bmatrix} 0 & 0 \\ 0 & 1 \end{bmatrix}\right) \in \Sigma_1, \;\; \Leftrightarrow \;\; a_1 < -0.5, a_2 \in \mathbb{R}.
$$

The Lyapunov inequality with a diagonal certificate reads as 
    \begin{align} \label{E:ExampleLyaDiag}
    &\begin{bmatrix} p_1 & 0\\ 0 & p_2 \end{bmatrix} \begin{bmatrix} 1 & 2\\ a_1 & a_2-k_2 \end{bmatrix}  + \begin{bmatrix} 1 & 2\\ a_1 & a_2-k_2 \end{bmatrix}^T \begin{bmatrix} p_1 & 0\\ 0 & p_2 \end{bmatrix}  \nonumber \\
    = & \begin{bmatrix} 2p_1 & 2p_1+a_1p_2\\ 2p_1 + a_1p_2 & 2a_2p_2 - 2k_2p_2 \end{bmatrix} \prec 0,
    \end{align}
where $p_1 > 0, p_2 > 0$. Since the first principle minor $2p_1 > 0$, we know that~\eqref{E:ExampleLyaDiag} is infeasible, $\forall a_1, a_2, k_2$.
Thus, we have
$$
    \left( \begin{bmatrix} 1 & 2\\ a_1 & a_2 \end{bmatrix}, \begin{bmatrix} 0 & 0 \\ 0 & 1 \end{bmatrix}\right) \notin \Sigma_2, \forall a_1 \in \mathbb{R}, a_2 \in \mathbb{R}.
$$

If both subsystems are fully actuated, \emph{i.e.}, $B_{1} = 1, B_{2} = 1$ in~\eqref{E:Example}, then according to Proposition~\ref{proposition:fullyActuated}, we know
$$
    \begin{aligned}
    \left( \begin{bmatrix} 1 & 2\\ a_1 & a_2 \end{bmatrix}, \begin{bmatrix} 1 & 0 \\ 0 & 1 \end{bmatrix}\right) &\in \Sigma_0, \forall a_1 \in \mathbb{R}, a_2 \in \mathbb{R}. \\
    \left( \begin{bmatrix} 1 & 2\\ a_1 & a_2 \end{bmatrix}, \begin{bmatrix} 1 & 0 \\ 0 & 1 \end{bmatrix}\right) &\in \Sigma_1, \forall a_1 \in \mathbb{R}, a_2 \in \mathbb{R}. \\
    \left( \begin{bmatrix} 1 & 2\\ a_1 & a_2 \end{bmatrix}, \begin{bmatrix} 1 & 0 \\ 0 & 1 \end{bmatrix}\right) &\in \Sigma_2, \forall a_1 \in \mathbb{R}, a_2 \in \mathbb{R}.
    \end{aligned}
$$%
This simple example also shows that the ability of actuating the nodes 
is important for strongly decentralized stabilization.

\bibliographystyle{IEEEtran}
\bibliography{Reference}

\begin{thebibliography}{10}
\providecommand{\url}[1]{#1}
\csname url@samestyle\endcsname
\providecommand{\newblock}{\relax}
\providecommand{\bibinfo}[2]{#2}
\providecommand{\BIBentrySTDinterwordspacing}{\spaceskip=0pt\relax}
\providecommand{\BIBentryALTinterwordstretchfactor}{4}
\providecommand{\BIBentryALTinterwordspacing}{\spaceskip=\fontdimen2\font plus
\BIBentryALTinterwordstretchfactor\fontdimen3\font minus
  \fontdimen4\font\relax}
\providecommand{\BIBforeignlanguage}[2]{{%
\expandafter\ifx\csname l@#1\endcsname\relax
\typeout{** WARNING: IEEEtran.bst: No hyphenation pattern has been}%
\typeout{** loaded for the language `#1'. Using the pattern for}%
\typeout{** the default language instead.}%
\else
\language=\csname l@#1\endcsname
\fi
#2}}
\providecommand{\BIBdecl}{\relax}
\BIBdecl

\bibitem{siljak2011decentralized}
D.~D. Siljak, \emph{Decentralized Control of Complex Systems}.\hskip 1em plus
  0.5em minus 0.4em\relax Courier Corporation, 2011.

\bibitem{lunze1992feedback}
J.~Lunze, \emph{Feedback Control of Large-scale Systems}.\hskip 1em plus 0.5em
  minus 0.4em\relax Prentice Hall New York, 1992.

\bibitem{wang1973stabilization}
S.-H. Wang and E.~Davison, ``On the stabilization of decentralized control
  systems,'' \emph{IEEE Trans. Autom. Control}, vol.~18, no.~5, pp. 473--478,
  1973.

\bibitem{geromel1994decentralized}
J.~C. Geromel, J.~Bernussou, and P.~L.~D. Peres, ``Decentralized control
  through parameter space optimization,'' \emph{Automatica}, vol.~30, no.~10,
  pp. 1565--1578, 1994.

\bibitem{jovanovic2016controller}
M.~R. Jovanovi{\'c} and N.~K. Dhingra, ``Controller architectures: Tradeoffs
  between performance and structure,'' \emph{European Journal of Control},
  vol.~30, pp. 76--91, July,2016.

\bibitem{rotkowitz2006characterization}
M.~Rotkowitz and S.~Lall, ``A characterization of convex problems in
  decentralized control,'' \emph{IEEE Trans. Autom. Control}, vol.~51, no.~2,
  pp. 274--286, 2006.

\bibitem{lin2011augmented}
F.~Lin, M.~Fardad, and M.~R. Jovanovic, ``Augmented \protect{Lagrangian}
  approach to design of structured optimal state feedback gains,'' \emph{IEEE
  Trans. Autom. Control}, vol.~56, no.~12, pp. 2923--2929, 2011.

\bibitem{fazelnia2017convex}
G.~Fazelnia, R.~Madani, A.~Kalbat, and J.~Lavaei, ``Convex relaxation for
  optimal distributed control problems,'' \emph{IEEE Trans. Autom. Control},
  vol.~62, no.~1, pp. 206--221, 2017.

\bibitem{furieri2017value}
L.~Furieri and M.~Kamgarpour, ``The value of communication in designing robust
  distributed controllers,'' \emph{arXiv 1711.05324}, 2017.

\bibitem{langbort2010distributed}
C.~Langbort and J.-C. Delvenne, ``Distributed design methods for linear
  quadratic control and their limitations,'' \emph{IEEE Trans. Autom. Control},
  vol.~55, no.~9, pp. 2085--2093, 2010.

\bibitem{farokhi2013optimal}
F.~Farokhi, C.~Langbort, and K.~H. Johansson, ``Optimal structured static
  state-feedback control design with limited model information for
  fully-actuated systems,'' \emph{Automatica}, vol.~49, pp. 326--337, 2013.

\bibitem{shah2013optimal}
P.~Shah and P.~Parrilo, ``H2-optimal decentralized control over posets: A
  state-space solution for state-feedback,'' \emph{IEEE Trans. Autom. Control},
  vol.~58, no.~12, pp. 3084--3096, 2013.

\bibitem{ahmadi2018distributed}
M.~Ahmadi, M.~Cubuktepe, U.~Topcu, and T.~Tanaka, ``Distributed synthesis using
  accelerated {ADMM},'' in \emph{2018 Annual American Control Conference
  (ACC)}.\hskip 1em plus 0.5em minus 0.4em\relax IEEE, 2018, pp. 6206--6211.

\bibitem{wang2017separable}
Y.-S. Wang, N.~Matni, and J.~C. Doyle, ``Separable and localized system-level
  synthesis for large-scale systems,'' \emph{IEEE Transactions on Automatic
  Control}, vol.~63, no.~12, pp. 4234--4249, 2018.

\bibitem{agler1988positive}
J.~Agler, W.~Helton, S.~McCullough, and L.~Rodman, ``Positive semidefinite
  matrices with a given sparsity pattern,'' \emph{Linear Algebra. Appl.}, vol.
  107, pp. 101--149, 1988.

\bibitem{grone1984positive}
R.~Grone, C.~R. Johnson, E.~M. S{\'a}, and H.~Wolkowicz, ``Positive definite
  completions of partial hermitian matrices,'' \emph{Linear Algebra. Appl.},
  vol.~58, pp. 109--124, 1984.

\bibitem{fukuda2001exploiting}
M.~Fukuda, M.~Kojima, K.~Murota, and K.~Nakata, ``Exploiting sparsity in
  semidefinite programming via matrix completion \protect{I}: General
  framework,'' \emph{SIAM J. Optimiz.}, vol.~11, no.~3, pp. 647--674, 2001.

\bibitem{andersen2010implementation}
M.~S. Andersen, J.~Dahl, and L.~Vandenberghe, ``Implementation of nonsymmetric
  interior-point methods for linear optimization over sparse matrix cones,''
  \emph{Math. Prog. Computat.}, no. 3-4, pp. 167--201, 2010.

\bibitem{mason2014chordal}
R.~P. Mason and A.~Papachristodoulou, ``Chordal sparsity, decomposing {SDPs}
  and the {Lyapunov} equation,'' in \emph{American Control Conference (ACC),
  2014}.\hskip 1em plus 0.5em minus 0.4em\relax IEEE, 2014, pp. 531--537.

\bibitem{ZKSP2017scalable}
Y.~Zheng, M.~Kamgarpour, A.~Sootla, and A.~Papachristodoulou, ``Scalable
  analysis of linear networked systems via chordal decomposition,'' in
  \emph{2018 European Control Conference (ECC)}.\hskip 1em plus 0.5em minus
  0.4em\relax IEEE, 2018, pp. 2260--2265.

\bibitem{andersen2014robust}
M.~Andersen, S.~Pakazad, A.~Hansson, and A.~Rantzer, ``Robust stability
  analysis of sparsely interconnected uncertain systems,'' \emph{IEEE Trans.
  Autom. Control}, vol.~59, no.~8, pp. 2151--2156, 2014.

\bibitem{Zheng2017Scalable}
Y.~Zheng, R.~P. Mason, and A.~Papachristodoulou, ``Scalable design of
  structured controllers using chordal decomposition,'' \emph{IEEE Trans.
  Autom. Control}, vol.~63, no.~3, pp. 752--767, March 2018.

\bibitem{zheng2017convex}
Y.~Zheng, M.~Kamgarpour, A.~Sootla, and A.~Papachristodoulou, ``Convex design
  of structured controllers using block-diagonal {Lyapunov} functions,''
  \emph{Technical report, avialable at arXiv:1709.00695}, 2017.

\bibitem{tanaka2011bounded}
T.~Tanaka and C.~Langbort, ``The bounded real lemma for internally positive
  systems and h-infinity structured static state feedback,'' \emph{IEEE
  transactions on automatic control}, vol.~56, no.~9, pp. 2218--2223, 2011.

\bibitem{boyd1994linear}
S.~Boyd, L.~El~Ghaoui, E.~Feron, and V.~Balakrishnan, \emph{Linear Matrix
  Inequalities in System and Control Theory}.\hskip 1em plus 0.5em minus
  0.4em\relax Society for Industrial and Applied Mathematics, 1994.

\bibitem{alavian2014stabilizing}
A.~Alavian and M.~Rotkowitz, ``Stabilizing decentralized systems with arbitrary
  information structure,'' in \emph{Decision and Control (CDC), 2014 IEEE 53rd
  Annual Conference on}.\hskip 1em plus 0.5em minus 0.4em\relax IEEE, 2014, pp.
  4032--4038.

\bibitem{sootla2017block}
A.~Sootla, Y.~Zheng, and A.~Papachristodoulou, ``Block-diagonal solutions to
  {Lyapunov} inequalities and generalisations of diagonal dominance,'' in
  \emph{2017 IEEE 56th Annual Conference on Decision and Control (CDC)}.\hskip
  1em plus 0.5em minus 0.4em\relax IEEE, 2017, pp. 6561--6566.

\bibitem{carlson1992block}
D.~Carlson, D.~Hershkowitz, and D.~Shasha, ``Block diagonal semistability
  factors and lyapunov semistability of block triangular matrices,''
  \emph{Linear Algebra and its Applications}, vol. 172, pp. 1--25, 1992.

\bibitem{blair1993introduction}
J.~R. Blair and B.~Peyton, ``An introduction to chordal graphs and clique
  trees,'' in \emph{Graph theory and sparse matrix computation}.\hskip 1em plus
  0.5em minus 0.4em\relax Springer, 1993, pp. 1--29.

\bibitem{vandenberghe2014chordal}
L.~Vandenberghe and M.~S. Andersen, ``Chordal graphs and semidefinite
  optimization,'' \emph{Found. Trends Optim.}, vol.~1, pp. 241--433, 2014.

\bibitem{griewank1984existence}
A.~Griewank and P.~L. Toint, ``On the existence of convex decompositions of
  partially separable functions,'' \emph{Mathematical Programming}, vol.~28,
  no.~1, pp. 25--49, 1984.

\bibitem{boyd2011distributed}
S.~Boyd, N.~Parikh, E.~Chu, B.~Peleato, and J.~Eckstein, ``Distributed
  optimization and statistical learning via the alternating direction method of
  multipliers,'' \emph{Found. Trends Mach. Learn.}, vol.~3, pp. 1--122, 2011.

\bibitem{sturm1999using}
J.~F. Sturm, ``Using {SeDuMi} 1.02, a {MATLAB} toolbox for optimization over
  symmetric cones,'' \emph{Optim. Methods Softw.}, vol.~11, no. 1-4, pp.
  625--653, 1999.

\bibitem{lofberg2004yalmip}
J.~L{\"o}fberg, ``\protect{YALMIP}: A toolbox for modeling and optimization in
  matlab,'' in \emph{Proc. IEEE Int. Symp. Computer Aided Control Syst.
  Design}.\hskip 1em plus 0.5em minus 0.4em\relax IEEE, 2004, pp. 284--289.

\end{thebibliography}

\end{document}